\documentclass[11pt]{article}%[12pt,twoside,openright,openbib]{article}

\usepackage{latexsym,epsfig}
\usepackage{amssymb,amsmath,amsfonts}
\usepackage{exscale} %% bei Schriftgroessen 11pt, 12pt
\usepackage{ifthen}
\usepackage{longtable}
\usepackage{subfigure}
\usepackage{caption}
\usepackage{cite}
\usepackage{lscape}

\begin{document}
\bibliographystyle{gabialpha}
\protect\pagenumbering{arabic}
\setcounter{page}{1}
 
\newcommand{\Zt}{\rm}

\newcommand{\ba}{\begin{array}}
\newcommand{\ea}{\end{array}}
\newcommand{\pot}{{\cal P}}
\newcommand{\curv}{\cal C}
\newcommand{\ddt} {\mbox{$\frac{\partial  }{\partial t}$}}
\newcommand{\hl}{\sf}
\newcommand{\hd}{\sf}

\newcommand{\Ad}{\mbox{\rm Ad}}
\newcommand{\Adsm}{\mbox{{\rm \scriptsize Ad}}}
\newcommand{\ad}{\mbox{\rm ad}}
\newcommand{\adsm}{\mbox{{\rm \scriptsize ad}}}
\newcommand{\diag}{\mbox{\rm Diag}}
\newcommand{\sect}{\mbox{\rm sec}}
\newcommand{\id}{\mbox{\rm id}}
\newcommand{\idsm}{\mbox{{\rm \scriptsize id}}}
\newcommand{\eps}{\varepsilon}
\newcommand{\Summ}{P}

\newcommand{\aL}{\mathfrak{a}}
\newcommand{\bL}{\mathfrak{b}}
\newcommand{\mL}{\mathfrak{m}}
\newcommand{\kL}{\mathfrak{k}}
\newcommand{\gL}{\mathfrak{g}}
\newcommand{\nL}{\mathfrak{n}}
\newcommand{\hL}{\mathfrak{h}}
\newcommand{\pL}{\mathfrak{p}}
\newcommand{\uL}{\mathfrak{u}}
\newcommand{\lL}{\mathfrak{l}}

\newcommand{\kG}{{\tt k}}
\newcommand{\nG}{{\tt n}}

\newcommand{\Cart}{$G=K e^{\overline{\aL^+}} K$}
\newcommand{\Area}{\mbox{Area}}
\newcommand{\Hd}{\mbox{\rm Hd}}
\newcommand{\Hdim}{\mbox{\rm dim}_{\mbox{\rm \scriptsize Hd}}}
\newcommand{\Tr}{\mbox{\rm Tr}}
\newcommand{\bs}{{\cal B}}
\newcommand{\bv}{{\rm B}}

\newcommand{\nc}{{\cal N}}
\newcommand{\MM}{{\cal M}}
\newcommand{\Ch}{{\cal C}}
\newcommand{\clCh}{\overline{\cal C}}
\newcommand{\Sh}{\mbox{Sh}}
\newcommand{\smSh}{\mbox{{\rm \scriptsize Sh}}}
\newcommand{\Cnt}{\mbox{\rm C}}
\newcommand{\preim}{(\pi^F)^{-1}}

\newcommand{\NN}{\mathbb{N}} \newcommand{\ZZ}{\mathbb{Z}}
\newcommand{\QQ}{\mathbb{Q}} \newcommand{\RR}{\mathbb{R}}
\newcommand{\KK}{\mathbb{K}} \newcommand{\FF}{\mathbb{F}}
\newcommand{\CC}{\mathbb{C}} \newcommand{\EE}{\mathbb{E}}
\newcommand{\XX}{X}%{\mathfrak{X}}
\newcommand{\HH}{I\hspace{-2mm}H}
\newcommand{\norm}{\Vert\hspace{-0.35mm}|}
\newcommand{\Sph}{\mathbb{S}}
\newcommand{\ganz}{\overline{\XX}}
\newcommand{\rand}{\partial\XX}
\newcommand{\prodrand}{\partial\XX_1\times\partial\XX_2} 
\newcommand{\regrand}{\partial\XX^{reg}}
\newcommand{\singrand}{\partial\XX^{sing}}
\newcommand{\Frand}{\partial^F\XX}
\newcommand{\Lim}{L_\Gamma}          
\newcommand{\Flim}{F_\Gamma}
\newcommand{\reglim}{L_\Gamma^{reg}}
\newcommand{\radlim}{L_\Gamma^{rad}}
\newcommand{\raylim}{L_\Gamma^{ray}}
\newcommand{\horinf}{\mbox{Vis}^{\infty}}
\newcommand{\horF}{\mbox{Vis}^F}
\newcommand{\Sml}{\mbox{Small}}
\newcommand{\SmlF}{\mbox{Small}^F}

\newcommand{\ifl}{\qquad\Longleftrightarrow\qquad}
\newcommand{\at}{\!\cdot\!}
\newcommand{\ging}{\gamma\in\Gamma}
\newcommand{\xo}{{o}}
\newcommand{\gamo}{{\gamma\xo}}
\newcommand{\gam}{\gamma}
\newcommand{\gax}{h}
\newcommand{\gxi}{{G\!\cdot\!\xi}}
\newcommand{\bd}{$(b,\theta)$-densit}
\newcommand{\bt}{$(b,\theta)$-densit}
\newcommand{\cd}{$(\alpha,\Gamma\at\xi)$-density}
\newcommand{\be}{\begin{eqnarray*}}
\newcommand{\ee}{\end{eqnarray*}}

\newcommand{\an}{\ \mbox{and}\ }
\newcommand{\as}{\ \mbox{as}\ }
\newcommand{\diam}{\mbox{diam}}
\newcommand{\is}{\mbox{Is}}
\newcommand{\Ax}{\mbox{Ax}}
\newcommand{\Fix}{\mbox{Fix}}
\newcommand{\Par}{F}
\newcommand{\Min}{\mbox{Fix}}
\newcommand{\rel}{\mbox{Rel}_\Gamma}
\newcommand{\vol}{\mbox{vol}}
\newcommand{\Td}{\mbox{Td}}
\newcommand{\piF}{\pi^B}
\newcommand{\piKM}{\pi^I}

\newcommand{\for}{\ \mbox{for}\ }
\newcommand{\pr}{\mbox{pr}}
\newcommand{\sh}{\mbox{sh}}
\newcommand{\shi}{\mbox{sh}^{\infty}}
\newcommand{\rank}{\mbox{rank}}
\newcommand{\supp}{\mbox{supp}}
\newcommand{\mass}{\mbox{mass}}
\newcommand{\kernel}{\mbox{kernel}}
\newcommand{\st}{\mbox{such}\ \mbox{that}\ }
\newcommand{\Stab}{\mbox{Stab}}
\newcommand{\Root}{\Sigma}
\newcommand{\Cone}{\mbox{C}}
\newcommand{\wrt}{\mbox{with}\ \mbox{respect}\ \mbox{to}\ }
\newcommand{\where}{\ \mbox{where}\ }

\newcommand{\con}{{\sc Consequence}\newline}
\newcommand{\rem}{{\sc Remark}\newline}
\newcommand{\prf}{{\sl Proof.\  }}
\newcommand{\qed}{$\hfill\Box$}

\newenvironment{rmk} {\newline{\sc Remark.\ }}{}  
\newenvironment{rmke} {{\sc Remark.\ }}{}  
\newenvironment{rmks} {{\sc Remarks.\ }}{}  
\newenvironment{nt} {{\sc Notation}}{}  

\newtheorem{satz}{\bf Theorem}

\newtheorem{df}{\sc Definition}[section]
\newtheorem{cor}[df]{\sc Corollary}
\newtheorem{thr}[df]{\bf Theorem}
\newtheorem{lem}[df]{\sc Lemma}
\newtheorem{prp}[df]{\sc Proposition}
\newtheorem{ex}{\sc Example}
\newenvironment{pros}{{\sc Properties:}}
%\newtheorem{cn}[cor]{\sc Consequence}

%%% Local Variables: 
%%% mode: latex
%%% TeX-master: "Article"
%%% End: 

%\title{\sc Exponent of Growth of discrete isometry groups acting on a product of Hadamard manifolds}
\title{{\sc Generalized Patterson-Sullivan measures for products of Hadamard spaces}}%\\{\sc (Preliminary draft)}}
\author{\sc Gabriele Link\thanks{supported by the DFG grant LI 1701/1-1}}
\date{\today}
\maketitle
\begin{abstract} Let $\Gamma$ be a  discrete group acting by isometries  on a product  $\XX=\XX_1\times \XX_2$ of
Hadamard spaces. We further require that  $\XX_1$, $\XX_2$ are locally compact and %geodesically complete, and 
$\Gamma$ contains two  elements projecting to a pair of independent rank one isometries in each factor.  
Apart from discrete groups acting by isometries on a product of CAT$(-1)$-spaces, the probably most interesting 
examples of such groups are Kac-Moody groups over finite fields acting on the Davis complex of their associated twin building.  
%which are known to act as irreducible lattices on a product of Hadamard spaces, 
In \cite{MR2629900} we showed that the regular 
geometric limit set splits as a product $F_\Gamma\times P_\Gamma$, where $F_\Gamma\subseteq\rand_1\times \rand_2$ is the projection 
of the geometric limit set to $\rand_1\times \rand_2$, and $P_\Gamma$ encodes the ratios of the speed of divergence of orbit points in each factor.
%denotes the set of accumulation points of the ratios of 
%the distances of the projections of orbit points to the two factors. 
Our aim in this paper is a description of the limit set from a measure 
theoretical point of view. We first study the conformal density obtained from the classical Patterson-Sullivan construction, then generalize this
construction to obtain measures supported in each $\Gamma$-invariant subset of the regular limit set and investigate
their properties. Finally we show that the Hausdorff dimension of the radial limit set in each $\Gamma$-invariant subset of 
$\Lim$ is bounded above by the exponential growth rate introduced in \cite{MR2629900}.
\end{abstract}

\vspace{0.2cm}

\section{Introduction}

Let $(\XX_1,d_1)$, $(\XX_2,d_2)$ be Hadamard spaces, i.e. complete simply connected metric spaces of non-positive Alexandrov curvature, and $(\XX,d)$ 
the product $\XX_1\times \XX_2$ endowed with the metric $d=\sqrt{d_1^2+d_2^2}$. Assume moreover that $\XX_1$, $\XX_2$ are locally compact. Each metric
 space $\XX, \XX_1,\XX_2 $ can be compactified by adding its  geometric boundary $\rand$, $
\rand_1$, $\rand_2$   endowed with the cone topology (see \cite[Chapter II]{MR1377265}). It is well-known that the {\hl regular geometric boundary} 
$\regrand$ of $\XX$ -- which consists of the set of equivalence classes of geodesic rays which do not project to a point in one of the factors -- is a 
dense open subset of $\rand$ homeomorphic to $\rand_1\times \rand_2\times (0,\pi/2)$. The last factor in this product is called the {\hl slope} of a 
point in $\regrand$.  The singular geometric boundary $\singrand=\rand\setminus\regrand$ consists of two strata homeomorphic to $\rand_1$, $\rand_2$ 
respectively. We assign slope $0$ to the first and slope $\pi/2$ to the second one. 

For a group $\Gamma\subseteq \is(\XX_1)\times \is(\XX_2)$ acting properly discontinuously by isometries on $\XX$ the limit set is defined by 
$\Lim:=\overline{\Gamma\at x}\cap\rand$, where $x\in\XX$ is arbitrary. In order to relate the critical exponent of a Fuchsian group to the Hausdorff 
dimension of its limit set, S.~J.~Patterson (\cite{MR0450547})  and D.~Sullivan (\cite{MR556586}) developed a theory of conformal densities. It 
turned out that for higher rank symmetric spaces and Euclidean buildings these densities in general detect only a small part of the geometric limit 
set (see \cite{MR1675889}). In order to measure the limit set in each invariant subset of the limit set, a class of generalized conformal densities 
were independently introduced in  \cite{MR1935549} and \cite{MR2062761}.  One of the main goals in this paper is to adapt this construction to 
 discrete groups  $\Gamma\subseteq\is(\XX_1)\times\is(\XX_2)$ which contain a pair of isometries  projecting to independent rank one elements in each factor. 
Related questions were considered by M.~Burger (\cite{MR1230298}) for graphs of convex cocompact groups in a product of rank one symmetric spaces, and  
by F.~Dal'bo and I.~Kim (\cite{DalboKim}) for discrete isometry groups of a product of Hadamard manifolds of pinched negative curvature.  %
%However, in our  more general setting even the construction of a conformal density is new. 

One important class of examples satisfying our conditions are Kac-Moody groups  $\Gamma$ over a finite field which act by isometries on a product 
$\XX=\XX_1\times\XX_2$, the CAT$(0)$-realization of the associated twin building ${\cal B}_+\times{\cal B}_-$.  Indeed, there  exists an element 
$h=(h_1,h_2)$ projecting to a rank one element in each factor by  Remark 5.4 and the proof of Corollary~1.3 in \cite{MR2585575}. Moreover, the 
action of the Weyl group produces many  axial isometries $g=(g_1,g_2)$ with $g_i$ rank one and independent from $h_i$ for $i=1,2$. Notice that if 
the order of the ground field is sufficiently large, then $\Gamma\subseteq\is(\XX_1)\times\is(\XX_2)$ is an irreducible lattice (see e.g. 
\cite{MR1715140} and \cite{MR2485882}). 

A second type of examples are groups acting properly discontinuously on a product of locally compact Hadamard spaces of strictly negative Alexandrov 
curvature (compare \cite{DalboKim} in the manifold setting). In this special case every non-elliptic and non-parabolic isometry in one of the factors 
is already a rank one element.  Prominent examples here which are already covered by the above mentioned results of J.~F.~ Quint  and the author are 
Hilbert modular groups acting as irreducible lattices on a product of hyperbolic planes and graphs of convex cocompact groups of rank one symmetric 
spaces (see also \cite{MR1230298}). 

Before we can state our results we need some definitions. We fix a base point $\xo=(\xo_1,\xo_2)\in\XX$. For $\theta\in [0,\pi/2]$ we  denote 
$\rand_\theta$ the set of points in the geometric boundary of slope $\theta$. Moreover, for $i\in\{1,2\}$ and $\eta_i\in\rand_i$ let 
$\bs_{\eta_i}(\cdot, \xo_i)$ denote the Busemann function centered at $\eta_i$ based at $\xo_i$.

A central role throughout the paper is played by the exponent of growth of $\Gamma$ of given slope $\theta\in [0,\pi/2]$ introduced in \cite{MR2629900}. 
For  $n\in\NN$  and $\eps>0$  we put 
$$N_\theta^\eps(n):=\#\{\gamma=(\gamma_1,\gamma_2)\in\Gamma: n-1<d(\xo,\gamma\xo)\le n,\, \big|\arctan \frac{d_2(\xo_2,\gamma_2\xo_2)}{d_1(\xo_1,\gamma_1\xo_1)}
-\theta\big|<\eps\}\,.$$
\begin{df}\label{expgrowth}
The {\hd exponent of  growth} of $\Gamma$ of slope $\theta\in [0,\pi/2]$ is defined by
$$\delta_\theta(\Gamma):=\liminf_{\eps\to 0}\left( \limsup_{n\to\infty}\frac1{n}\log N_\theta^\eps(n)\right)\,.$$
\end{df}
The quantity $\delta_\theta(\Gamma)$ can be thought of as a function of $\theta\in [0,\pi/2]$ which describes the exponential growth rate of orbit points 
converging to limit points of slope $\theta$. It is an invariant of $\Gamma$ which carries more information than the critical exponent $\delta(\Gamma)$. 
Moreover, Theorem~7.4 in \cite{MR2629900} implies that there exists a unique slope $\theta_*\in [0,\pi/2]$ \st the exponent of growth 
of $\Gamma$ is maximal for this slope and equal to the critical exponent $\delta(\Gamma)$.  %$\delta_{\theta_*}(\Gamma)=\delta(\Gamma)$.

Our first result concerns the measures on the geometric boundary obtained by the classical Patterson-Sullivan construction. Similar to the case of higher rank 
symmetric spaces or Euclidean buildings we have the following result:\\[2mm]
 {\bf Theorem A}$\quad$ {\sl The Patterson-Sullivan construction produces a conformal density with support in a single $\Gamma$-invariant subset of the geometric 
limit set. Every point in its support has slope $\theta_*$ as above.} \\[-1mm]

Thus in order to measure the remaining $\Gamma$-invariant subsets of the limit set, we need a more sophisticated construction. Inspired by the paper \cite{MR1230298} of M.~Burger 
we  therefore consider densities with one more degree of freedom than the classical 
conformal density:  
\begin{df}\label{bdensi}
Let $\MM^+(\rand)$ denote the cone of positive finite Borel measures
on $\rand$, $\theta\in [0,\pi/2]$ and $b=(b_1,b_2)\in \RR^2$.
%We further require\newline
%$(*)\quad  b_1=0$ if $\theta=\pi/2$ and $b_2=0$ if $\theta=0$. 
A $\Gamma$-invariant {\hl $(b,\theta)$-density} is a
continuous map 
$$\begin{array}{ccc}\mu:\XX &\to &\MM^+(\rand)\\
x&\mapsto&\mu_x\end{array}$$ 
\st for any $x\in\XX$ the following three properties hold:%\\[1mm]
\begin{flushleft}
\begin{tabular}{rl}
 {\rm (i)}\ &   $\emptyset\neq\supp(\mu_{x})\subseteq\Lim\cap \rand_\theta$,\\[1mm]
{\rm (ii)}\  & $\forall \,\gamma\in\Gamma\qquad \gamma_*\mu_x=\mu_{\gamma  x}$\footnotemark[1],\\[1mm]
{\rm (iii)}\ & if $\ \theta\in (0,\pi/2),\ $ then $\quad \forall\, \tilde\eta=(\eta_1,\eta_2,\theta)\in\supp(\mu_\xo)$\\[1mm]
&  \hspace{5cm} $ \displaystyle{\frac{d\mu_x}{d\mu_\xo}}(\tilde\eta)=e^{b_1\bs_{\eta_1}(\xo_1,x_1)+b_2
\bs_{\eta_2}(\xo_2,x_2)}$,\\[2mm]
& if $\ \theta=0,\ $ then $\ \, b_2=0\ \an \ \, \forall\, \tilde\eta=\eta_1\in\supp(\mu_\xo)\quad \displaystyle{\frac{d\mu_x}{d\mu_\xo}}(\tilde\eta)
=e^{b_1\bs_{\eta_1}(\xo_1,x_1)}$,\\[3mm]
& if $\ \theta=\frac{\pi}2,\ $ then $\ \, b_1=0\ \an \ \, \forall\, \tilde\eta=\eta_2\in\supp(\mu_\xo)
\quad \displaystyle{\frac{d\mu_x}{d\mu_\xo}}(\tilde\eta)=e^{b_2\bs_{\eta_2}(\xo_2,x_2)}$.
\end{tabular}
\footnotetext[1]{Here $\gamma_*\mu_x$ denotes the measure defined by $\gamma_*\mu_x(E)=\mu_x(\gamma^{-1}E)$ for any Borel set $E\subseteq\rand$}
\end{flushleft}
\end{df} 

Notice that the conformal density from Theorem A is a special case of such a density with support in $\rand_{\theta_*}$ and parameters $b_1=\delta(\Gamma)\cos\theta_*$, $b_2=\delta(\Gamma)\sin\theta_*$. %So all the results concerning $(b,\theta)$-densities in particular apply to any conformal density supported in a single $\Gamma$-invariant subset of the geometric limit set.

We next give a criterion for the existence of a \bd y.\\[2mm]
{\bf Theorem B}$\quad$ {\sl If $\theta\in (0,\pi/2)$ is such that $\delta_\theta(\Gamma)>0$, then there exists a \bd y for some parameters 
$b=(b_1,b_2)\in\RR^2$.}\\[-1mm]

In Section~\ref{genPseries} we will explicitly describe the construction of such a  \bd y. Notice that our method does not cover the cases $\theta=0$ and $\theta=\pi/2$ in general.
 However, if $\delta_\theta(\Gamma)=\delta(\Gamma)$, then by Theorem A the classical Patterson-Sullivan construction provides a \bd y, whether $\theta$ belongs to $(0,\pi/2)$ or not. 
Unfortunately, we do not know of an example with $\delta_\theta(\Gamma)=\delta(\Gamma)$ for $\theta=0$ or $\theta=\pi/2$.

The following results  about $(b,\theta)$-densities in particular apply to any conformal density supported in a single $\Gamma$-invariant subset of the geometric limit set, not 
only the one obtained by the classical Patterson-Sullivan construction. Our 
main tool is a so-called shadow lemma for $(b,\theta)$-densities, which is a generalization of the well-known shadow lemma for conformal densities. It first gives a condition 
for the parameters of a $(b,\theta)$-density in terms of the exponent of growth. \\[2mm]
{\bf Theorem C}$\quad$ {\sl If a $\Gamma$-invariant $(b,\theta)$-density exists for some  $\theta\in (0,\pi/2)$, then \\[-3mm]

$\hspace{4cm}\delta_\theta(\Gamma)\le b_1\cos\theta +b_2\sin\theta\,.$}\\[-1mm]
%Condition $(*)$ is necessary if $\xi\in\singrand$ in order for (iii) to be well defined

The following subsets of the geometric limit set will play an important role in the sequel.
\begin{df}\label{raddef}
A point $\tilde \xi\in \rand$  is called a {\hd radial limit point} of $\Gamma$ if there exists a sequence $(\gamma_n)=\big((\gamma_{n,1},\gamma_{n,2})\big)\subset\Gamma$ 
\st $\gamma_n\xo$ converges to $\tilde\xi$, and the following condition holds:

If $\tilde\xi=(\xi_1,\xi_2,\theta)\in \regrand$, then  for $i\in\{1,2\}$ $\gamma_{n,i}\xo_i$ stays at bounded distance of one (and hence any) 
geodesic ray in the class of $\xi_i$, if $i=1,2$ and  $\tilde \xi=\xi_i\in \rand_i\subset\singrand$, then $\gamma_{n,i}\xo_i$ stays at bounded distance of one (and hence any) geodesic ray in the class of $\xi_i$.

 We will denote the set of radial limit points of $\Gamma$ by $\radlim$. % In the more restrictive case that $\gamma_n\xo$ stays at bounded distance of a geodesic ray  in
% the class of $\tilde\xi$ in the product $\XX$, we call $\tilde\xi$ a {\hd ray limit point}.
\end{df}
Notice that in general, a radial limit point $\tilde \xi$ is not  approached by a sequence $\gamma_n\xo$, $\gamma_n\in\Gamma$, at bounded distance of a geodesic ray in the class of 
$\tilde \xi$.   % this condition is weaker than the condition that $\gamma_n\xo$ stays at bounded distance of a geodesic ray in the class of $\xi$.   

Our next statement shows that for certain \bd ies the corresponding exponent of growth $\delta_\theta(\Gamma)$ is completely determined by the parameters $\theta\in (0,\pi/2)$ 
and $b=(b_1,b_2)\in\RR^2$.\\[2mm] %are closely related to the exponent of growth.\\[2mm]
{\bf Theorem D}$\quad$ {\sl If $\theta\in (0,\pi/2)$ and $\mu$ is a $\Gamma$-invariant  \bd y which gives positive measure to the radial limit set, then 
$\ \delta_\theta(\Gamma)=b_1\cos\theta +b_2\sin\theta$.}\\[3mm]
%\begin{satz}
%If $\theta\in (0,\pi/2)$ and $\mu$ is a $\Gamma$-invariant  \bd y, then a radial limit point is not a point mass for $\mu$.
%\end{satz}
The following theorem gives a restriction for  the atomic part of our measures.\\[2mm] 
{\bf Theorem E}$\quad$ {\sl If $\theta\in (0,\pi/2)$ \st  $\delta_\theta(\Gamma)>0$, and $\mu$ is a $\Gamma$-invariant  \bd y, then a radial limit point is not a point mass for $\mu$. 
}\\[-1mm]

Moreover, using a Hausdorff measure on the geometric boundary as proposed by G.~Knieper (\cite[Section 4]{MR1465601}), we have the following\\[2mm] 
{\bf Theorem F}$\quad$ {\sl For any $\theta\in [0,\pi/2]$ with $\Lim\cap\rand_\theta\ne\emptyset$ we have \\[-4mm]

$\hspace{4cm} \Hdim(\radlim\cap \rand_\theta)\le\delta_{\theta}(\Gamma)\,.$ }\\[-1mm]

Unfortunately, a precise estimate for the Hausdorff dimension can be given only for a particular class of groups which we choose to 
call {\hl radially cocompact}. Examples of such groups are uniform lattices %cocompact groups 
and products of convex cocompact groups acting on a product of Hadamard manifolds of pinched negative curvature. Since Kac-Moody groups over finite fields are never cocompact, 
we do not know whether a similar result holds for them.\\[2mm]
{\bf Theorem G}$\quad$ {\sl If $\ \Gamma$ is  radially cocompact and
$\theta\in (0,\pi/2)$ such that $\delta_\theta(\Gamma)>0$, then \\[-4mm]

$\hspace{4cm}\Hdim(\radlim\cap \rand_\theta)=\delta_{\theta}(\Gamma)\,.$}\\[-1mm]

The paper is organized as follows: 
In Section~\ref{Prelim} we recall  basic facts about Hadamard spaces and rank one isometries. Section~\ref{prodHadspaces} deals with the product case and provides some tools for the proof 
of the so-called shadow lemma in Section 7.  In Section~\ref{ExpGrowth} we introduce and study the properties of the exponent of growth. Section~\ref{ClasPatSulconst} recalls the 
classical Patterson-Sullivan construction in our setting. The main new result here is  Theorem A. In Section~\ref{genPseries} we introduce a generalized Poincar{\'e} series that  
allows to construct \bd ies, and therefore proves Theorem B. Using the shadow lemma, in Section~\ref{Propbdies} we deduce properties of \bd ies and prove Theorems C, D and E. Section~\ref{Hausdorff} 
finally is concerned with  the Hausdorff dimension of the limit set and the proofs of Theorems F and G.\\[1mm]

{\bf Acknowledgements:}  This paper was initiated during the author's stay at IHES in Bures-sur-Yvette. She warmly thanks the institute for its hospitality and the inspiring atmosphere.

\section{Preliminaries}\label{Prelim}

The purpose of this section is to introduce terminology and notation and to summarize basic results about Hadamard spaces and rank one isometries. %y groups of such. 
The main references here  are \cite{MR1744486} and  \cite{MR1377265} (see also \cite{MR1383216}, and \cite{MR823981},\cite{MR656659} in the case of Hadamard manifolds). 

Let $(\XX,d)$ be a metric space. A {\hl geodesic path} joining $x\in\XX$ to $y\in\XX$  is a map $\sigma$ from a closed interval $[0,l]\subset \RR$ to $\XX$ \st $\sigma(0)=x$, 
$\sigma(l)=y$ and $d(\sigma(t), \sigma(t'))=|t-t'|$ for all $t,t'\in [0,l]$.  We will denote such a geodesic path $\sigma_{x,y}$. $\XX$ is called {\hl geodesic} if any two points 
in $\XX$ can be connected by a geodesic path, if this path is unique we say that $\XX$ is {\hl uniquely geodesic}. In this text $\XX$ will be a Hadamard space, i.e. a complete 
geodesic metric space in which all triangles satisfy the CAT$(0)$-inequality. This implies in particular that $\XX$ is simply connected and uniquely geodesic.   A {\hl geodesic} 
or {\hl geodesic line} in $\XX$ is a map $\sigma:\RR\to\XX$ \st $d(\sigma(t), \sigma(t'))=|t-t'|$ for all $t,t'\in\RR$, a {\hl geodesic ray} is a map $\sigma:[0,\infty)\to \XX$ 
\st $d(\sigma(t), \sigma(t'))=|t-t'|$ for all $t,t'\in [0,\infty)$. Notice that in the non-Riemannian setting completeness of $\XX$ does not imply geodesically completeness, i.e. 
not every geodesic path or ray can be extended to a geodesic.

%$\XX$ is called {\hl geodesically complete}, if everyÄ geodesic path can be extended to a map $\sigma:\RR\to \XX$ \st $d(\sigma(t), \sigma(t'))=|t-t'|$ for all $t,t'\in\RR$. 
%Such a map will be called a {\hl geodesic} or {\hl geodesic line} in $\XX$. A {\hl geodesic ray} is a map $\sigma:[0,\infty)\to \XX$ \st $d(\sigma(t), \sigma(t'))=|t-t'|$ 
%for all $t,t'\in [0,\infty)$.

From here on we will assume that $\XX$ is a locally compact Hadamard space.  The geometric boundary $\rand$ of
$\XX$ is the set of equivalence classes of asymptotic geodesic rays endowed with the cone topology (see e.g. \cite[chapter~II]{MR1377265}). The action of the isometry group 
$\is(\XX)$ on $\XX$ naturally extends to an action by homeomorphisms on the geometric boundary. Moreover, since $\XX$ is locally compact, this boundary $\rand$ is compact 
and the space $\XX$ is a dense and open subset of the compact space $\ganz:=\XX\cup\rand$.  For $x\in\XX$ and $\xi\in\rand$ arbitrary there exists a  geodesic ray emanating 
from $x$ which belongs to the class of $\xi$. We will denote such a ray $\sigma_{x,\xi}$.

We say that two points $\xi$, $\eta\in\rand$ can be joined by a geodesic if there exists a geodesic $\sigma:\RR\to\XX$ \st $\sigma(-\infty)=\xi$ and 
$\sigma(\infty)=\eta$. It is well-known that if $\XX$ is  CAT$(-1)$, i.e. of negative Alexandrov curvature bounded above by $-1$, then every pair of
distinct points in the geometric boundary can be joined by a geodesic. This is not true in the general CAT$(0)$-case. 

Let $x, y\in \XX$, $\xi\in\rand$ and $\sigma$ a geodesic ray in the class of $\xi$. We put 
\begin{equation}\label{buseman}
 \bs_{\xi}(x, y)\,:= \lim_{s\to\infty}\big(
d(x,\sigma(s))-d(y,\sigma(s))\big)\,.
\end{equation}
This number is independent of the chosen ray $\sigma$, and the function
\be \bs_{\xi}(\cdot , y):
\quad \XX &\to & \RR\\
x &\mapsto & \bs_{\xi}(x, y)\ee
is called the {\hl Busemann function} centered at $\xi$ based at $y$ (see also \cite{MR1377265}, chapter~II). From the definition one immediately obtains the following 
properties of the Busemann function:
\begin{eqnarray}
\bs_{\xi}(x,y)&=& -\bs_\xi(y,x)\\
|\bs_{\xi}(x, y)|&\le &d(x,y)\label{bounded}\\
\bs_{\xi}(x, z)&=&\bs_{\xi}(x, y)+\bs_{\xi}(y,z)\label{cocycle}\\
\bs_{g\cdot\xi}(g\at x,g\at y) & =& \bs_{\xi}(x, y)\nonumber
\end{eqnarray}
for all $x,y,z\in\XX$, $\xi\in\rand$ and $g\in\is(\XX)$. Moreover, $\bs_{\xi}(x,y)=d(x,y)$ if and only if $y$ is a point on the geodesic ray $\sigma_{x,\xi}$, and we 
have the following easy 
\begin{lem}\label{distbusbounded}
Let $c>0$, $x,y\in \XX$, $x\ne z$, and $\xi\in\rand$ \st $d(z,\sigma_{x,\xi})<c$. Then 
$$0\le d(x,z)-\bs_\xi(x,z)<2c\,.$$
\end{lem}
\prf\  The first inequality is (\ref{bounded}). For the second one let $y\in\XX$ be a point on the geodesic ray $\sigma_{x,\xi}$ \st $d(z,y)<c$. Then for all $s>d(x,y)$ 
we have by the triangle inequality
\be d(x,\sigma_{x,\xi}(s))-d(z,\sigma_{x,\xi}(s))& \ge & d(x,\sigma_{x,\xi}(s))-d(z,y)-d(y,\sigma_{x,\xi}(s))\\
&=& d(x,y)-d(z,y)> d(x,y)-c\,,
\ee
hence $\ d(x,z)-\bs_\xi(x,z)\le d(x,y)+c -d(x,y)+ c=2c$.\qed\\[-1mm]
%The notion of rank one geodesics and rank one isometries is by now standard, but we recall the definitions for the convenience of the reader. 

A geodesic $\sigma: \RR\to\XX$  is said to {\hl bound a flat half-plane}  if there exists a closed convex subset $i([0,\infty)\times \RR)$ in $\XX$ isometric to 
$[0, \infty)\times\RR$ \st $\sigma(t)=i(0,t)$ for all $t\in \RR$. Similarly, a geodesic $\sigma: \RR\to\XX$  bounds a {\hl flat strip} of width $c>0$   if there 
exists a closed convex subset $i([0,c]\times \RR)$ in $\XX$ isometric to $[0, c]\times\RR$ \st $\sigma(t)=i(0,t)$ for all $t\in \RR$.
% \marginpar{defined in singular case???}{totally geodesic}  isometric embedding $i:[0, c]\times\RR\to\XX$ such that $i(0, t) =\sigma(t)$ for any $t\in\RR$.
We call a geodesic $\sigma:\RR\to\XX$ a {\hl rank one geodesic} if  $\sigma$ does not bound a
flat half-plane. 

The following important lemma states that even though we cannot join any two distinct points in the geometric boundary of $\XX$, given a rank one geodesic we can 
at least join points in a neighborhood of its extremities. More precisely, we have the following well-known
\begin{lem}\label{joinrankone} (\cite{MR1377265}, Lemma III.3.1)\ 
Let $\sigma:\RR\to\XX$ be a rank one geodesic. Then there exist $c>0$ and neighborhoods $U$ of $\sigma(-\infty)$ and $V$ of $\sigma(\infty)$ in 
$\ganz$ \st for any $\xi\in U$ and $\eta \in V$ there exists a rank one geodesic joining $\xi$ and $\eta$. Moreover, any such geodesic $\sigma'$ satisfies 
$d(\sigma', \sigma(0))\le c$..
\end{lem}

The following kind of isometries will
 play a central role in the sequel. 
\begin{df} \label{axialisos}
An isometry $h$ of $\XX$ is called {\hd axial}, if there exists a constant
$l=l(h)>0$ and a geodesic $\sigma$ \st
$h(\sigma(t))=\sigma(t+l)$ for all $t\in\RR$. We call
$l(h)$ the {\hd translation length} of $h$, and $\sigma$
an {\hd axis} of $h$. The boundary point
$h^+:=\sigma(\infty)$ is called the {\hd attractive fixed
point}, and $h^-:=\sigma(-\infty)$ the {\hd repulsive fixed
point} of $h$. We further put
$\Ax(h):=\{ x\in\XX\;|\, d(x,h x)=l(h)\}$.
\end{df}

We remark that $\Ax(h)$ consists of the union of parallel geodesics
translated by $h$, and 
$\overline{\Ax(h)}\cap\rand$ is exactly the set of fixed points of
$h$.
\begin{df}
An axial isometry  is called {\hd rank one} if it possesses a rank one axis. Two rank one isometries are called 
{\hd independent}, if their fixed point sets are disjoint. 
\end{df}

Notice that if $h$ is rank one, then $h^+$ and $h^-$ are the only fixed points of $h$. 
Let us recall the north-south dynamics of rank one isometries.
\begin{lem}\label{dynrankone}(\cite{MR1377265}, Lemma III.3.3)\ 
Let $h$ be a rank one isometry. Then
\begin{enumerate}
\item[(a)] Any $\xi\in\rand\setminus\{h^+\}$ can be joined to $h^+$ by a geodesic, and every geodesic joining  $\xi$ to $h^+$ is rank one, 
\item[(b)] given neighborhoods $U$ of $h^-$ and $V$ of $h^+$ in $\ganz$ 
there exists $N_0\in\NN$ \st\\
 $h^{-n}(\ganz\setminus V)\subset U$ and
$h^{n}(\ganz\setminus U)\subset V$ for all $n\ge N_0$.
\end{enumerate}
\end{lem}

If  $\Gamma$ is a group acting %properly discontinuously 
by isometries on a locally compact Hadamard space $\XX$ we define its {\hl geometric limit set}  by
$\Lim:=\overline{\Gamma\at x}\cap\rand$, where $x\in \XX$ is arbitrary.

\section{Products of Hadamard spaces}\label{prodHadspaces}

Now let $(\XX_1,d_1)$, $(\XX_2,d_2)$ be locally compact Hadamard spaces,  and $\XX=\XX_1\times \XX_2$  the product space endowed with the product 
distance $d=\sqrt{d_1^2+d_2^2}$. Notice that such a product is again a locally compact Hadamard space.  To any pair of points $x=(x_1,x_2)$, $z=(z_1,z_2)\in \XX$ 
we associate the vector 
\begin{equation}\label{distancevector}
H(x,z):= \left(\begin{array}{c} d_1(x_1,z_1)\\ d_2(x_2,z_2)\end{array}\right)\in \RR^2\,,
\end{equation} 
which we call the {\hl distance vector} of the pair $(x,z)$. If $z\neq x$ we further define the {\hl direction} of $z$  \wrt $x$ by 
\begin{equation}\label{slopeangle}
\theta(x,z):=\arctan\frac{d_2(x_2,z_2)}{d_1(x_1,z_1)}\,.
\end{equation}
Moreover, for convenience we set $\theta(x,x)=0$ for $x\in\XX$.

Clearly we have $H(z,x)=H(x,z)$ and $\theta(z,x)=\theta(x,z)$. Notice that we can also write 
$$H(x,z) = d(x, z)\left(\begin{array}{c} \cos \theta(x,z)\\ \sin\theta(x,z)\end{array}\right)\,,$$ 
hence in particular $\Vert H(x,z)\Vert =d(x,z)$, where $\Vert \cdot\Vert$ denotes the Euclidean norm in $\RR^2$. The following easy lemma states that distance
 vectors and directions are invariant by $\is(\XX_1)\times\is(\XX_2)$.
\begin{lem}\label{Ginvariance}
If $g=(g_1,g_2)\in \is(\XX_1)\times\is(\XX_2)$, $x=(x_1,x_2)$, $z=(z_1,z_2)\in\XX$, then 
$$ H(gx,gz)=H(x,z)\quad\mbox{and}\qquad \theta(gx,gz)=\theta(x,z)\,.$$
\end{lem}
\prf\ Since $g_1\in\is(\XX_1)$ and $g_2\in\is(\XX_2)$ we have 
$$ d_1(g_1x_1,g_1z_1)=d_1(x_1,z_1)\quad\mbox{and}\qquad d_2(g_2x_2,g_2z_2)=d_2(x_2,z_2)\,.$$ 
Hence by (\ref{distancevector}) and (\ref{slopeangle})
\be H(gx,gz) &= & \left(\begin{array}{c} d_1(g_1 x_1,g_1z_1)\\ d_2(g_2x_2,g_2z_2)\end{array}\right)= \left(\begin{array}{c} d_1(x_1,z_1)\\ d_2(x_2,z_2)\end{array}\right)=H(x,z)\,,\\
\theta(gx,gz)&=&\arctan\frac{d_2(g_2x_2,g_2z_2)}{d_1(g_1x_1,g_1z_1)}=\arctan\frac{d_2(x_2,z_2)}{d_1(x_1,z_1)}=\theta(x,z)\,.\hspace{28mm}\Box\ee

Denote $p_i:\XX\to \XX_i$, $i=1,2$, the natural projections. Every %unit speed 
geodesic path $\sigma:[0,l]\to\XX\,$ can be written as a product $\sigma(t)=(\sigma_1(t \cos \theta), \sigma_2(t \sin\theta))$, where $\theta\in [0,\pi/2]$ and 
$\sigma_1:[0,l \cos\theta]\to\XX_1$, $\sigma_2:[0,l \sin\theta]\to\XX_2$ are geodesic paths in $\XX_1$, $\XX_2$. $\theta$ equals the direction of $\sigma(l)$ with 
respect to $\sigma(0)$ and  is called the {\hl slope of $\sigma$}.  We say that a geodesic path $\sigma$ is  {\hl regular} if its slope is contained in the open 
interval $(0,\pi/2)$. In other words, $\sigma$ is regular if neither $p_1(\sigma([0,l]))$ nor $p_2(\sigma([0,l]))$ is a point.

If  $x\in\XX$  and $\sigma:[0,\infty)\to\XX$ is an arbitrary geodesic ray, then by elementary geometric estimates one has the relation 
\begin{equation}\label{RelSlopeDir}
 \theta = \lim_{t\to\infty} \theta(x,\sigma(t)) 
 \end{equation}
 between the slope $\theta$ of $\sigma$ and the directions of $\sigma(t)$, $t>0$, \wrt $x$. 
 Similarly, one can easily show that any two geodesic rays representing the same (possibly singular) point in the geometric boundary necessarily have the same slope. 
So  we may define the slope $\theta(\tilde\xi)$ of a point $\tilde\xi\in\rand$ as the slope of an arbitrary geodesic ray representing $\tilde\xi$. 
 
It is easy to see that a pair of distinct boundary points cannot be joined by a geodesic if they do not have the same slope. Moreover, two regular geodesic rays  $\sigma$, $\sigma'$  
of the same slope represent the same point in the geometric boundary if and only if  $\sigma_1(\infty)=\sigma_1'(\infty)$ and $\sigma_2(\infty)=\sigma_2'(\infty)$.  The regular geometric 
boundary $\regrand$ of $\XX$ is defined as the set of equivalence classes of regular geodesic rays and hence  is  homeomorphic  to $\rand_1\times \rand_2\times (0,\pi/2)$. 

 If $\gamma\in\is(\XX_1)\times \is(\XX_2)$ , then   the slope of $\gamma\at\tilde\xi$ equals the slope of $\tilde\xi$. In other words, if $\rand_\theta$ denotes the set of points in the 
geometric boundary of slope $\theta\in [0,\pi/2]$, then $\rand_\theta$ is invariant by the action of $\is(\XX_1)\times\is(\XX_2)$. Notice that points in $\singrand:= (\rand)_0\sqcup(\rand)_{\pi/2}$ 
are equivalence classes of geodesic rays which project to a point in one of the factors of $\XX$. Hence $(\rand)_0$ is homeomorphic to $\partial \XX_1$ and $(\rand)_{\pi/2}$ is homeomorphic to 
$\rand_2$. If  $\theta\in (0,\pi/2)$, then the set $\rand_\theta\subset\regrand$ is homeomorphic to the product $\rand_1\times\rand_2$. 
In the sequel we will often use the identification $\rand=\rand_1\sqcup \rand_2\sqcup \big(\rand_1\times \rand_2 \times (0,\pi/2)\big)$.

%We remark that a necessary condition for a sequence $(y_n)=\big((y_{n,1},y_{n,2})\big)\subset\XX$ to converge to a boundary point $\tilde\eta\in\rand$ is $\theta(\xo,y_n)\to\theta$ as $n\to\infty$. 
We remark that if $\theta\in (0,\pi/2)$, then a sequence $(y_n)=\big((y_{n,1},y_{n,2})\big)\subset\XX$ converges to $\tilde\eta=(\eta_1,\eta_2,\theta)$ 
if and only if $y_{n,1}\to \eta_1$, $y_{n,2}\to\eta_2$ and $\theta(\xo,y_n)\to\theta$ as $n\to\infty$. Similarly, 
$(y_n)=\big((y_{n,1},y_{n,2})\big)\subset\XX$ converges to 
$\tilde\eta=\eta_1\in (\rand)_0\cong\rand_1$ if and only if $y_{n,1}\to \eta_1$ and $\theta(\xo,y_n)\to 0$ as $n\to\infty$, 
and $(y_n)=\big((y_{n,1},y_{n,2})\big)\subset\XX$ converges to 
$\tilde\eta=\eta_2\in (\rand)_{\pi/2}\cong\rand_2$ if and only if $y_{n,2}\to \eta_2$ and $\theta(\xo,y_n)\to \pi/2$ as $n\to\infty$.

For higher rank symmetric spaces and Bruhat-Tits buildings there is a well-known notion of Furstenberg boundary which -- for a product of rank one spaces --
coincides with the product of the geometric boundaries. In our setting we choose to call the product $\rand_1\times \rand_2$ endowed with
the product topology the {\hl Furstenberg boundary} $\Frand$ of $\XX$. 
Using the above parametrization of $\regrand$ we have a natural projection
$$ \hspace{1cm} \ba{rcl}\pi^F\,:\qquad\ \regrand&\to& \Frand\\
(\xi_1,\xi_2,\theta)&\mapsto & (\xi_1,\xi_2)\, \ea $$
and a natural action of the  group $\is(\XX_1)\times\is(\XX_2)$ by homeomorphisms on the Furstenberg boundary of $\XX=\XX_1\times\XX_2$.

We say that two points $\xi=(\xi_1,\xi_2)$, $\eta=(\eta_1,\eta_2)\in\Frand$ are {\hl opposite} if $\xi_1$ and $\eta_1$ can be joined by a geodesic in $\XX_1$, and $\xi_2$, $\eta_2$ can be 
joined by a geodesic in $\XX_2$. Notice that if $\xi=(\xi_1,\xi_2)$, $\eta=(\eta_1,\eta_2)\in\Frand$ are opposite, then for any $\theta\in (0,\pi/2)$ the pair of points 
$\tilde\xi= (\xi_1,\xi_2,\theta)$, $\tilde\eta=(\eta_1,\eta_2,\theta)\in\rand_\theta$ provides a pair of boundary points which can be joined by a geodesic in $\XX$. The same holds for the 
pairs $\tilde\xi=\xi_1$, $\tilde\eta=\eta_1\in (\rand)_0$ and $\tilde\xi=\xi_2$, $\tilde\eta=\eta_2\in (\rand)_{\pi/2}$. 
%In the particular case that both  $\XX_1$, $\XX_2$ are CAT$(-1)$, this holds if and only if $\xi_1\neq \eta_1$ and $\xi_2\neq \eta_2$.

Let $\sigma_1$, $\sigma_2$ be geodesics in $\XX_1$,  $\XX_2$. We will call a  set $F$ in $\XX$ of the form $$F=\{(\sigma_1(t_1),\sigma_2(t_2)): t_1,t_2\in\RR\}\,$$
a  {\hl flat} in $\XX$. Notice that a flat defined in this way is a particular case of a geometric $2$-flat in $\XX$, i.e. a closed convex subset of $\XX$ isometric to Euclidean $2$-space. 
If $\XX_1$ and $\XX_2$ are CAT$(-1)$, then every geometric 2-flat in $\XX$ is a flat according to our definition. In general, however, this is not the case, because there may exist  
geometric $2$-flats in each factor. In particular, $\XX$ may contain geometric $n$-flats of larger dimensions $n$. 

The boundary $\partial F$ of a flat  $F$ is %\marginpar{only "contains"?}
topologically a circle, and $\pi^F(\partial F\cap\regrand)$ consists of 4 points. We say that a flat $F$ {\hl joins} $\xi$, $\eta\in\Frand$, $\xi\neq \eta$, 
if $\xi,\eta\in\pi^F(\partial F\cap \regrand)$. Notice that even if $\xi$, $\eta\in\Frand$ can be joined by a flat, $\xi$ and $\eta$ need not be opposite. 

It is easy to see that if $\XX_1$, $\XX_2$ are CAT$(-1)$, then any two distinct points $\xi$, $\eta$ in the Furstenberg boundary can be joined by a unique 
flat, analogously to the situation in higher rank symmetric spaces. This is clearly not true in general.

For $x=(x_1,x_2)$, $z=(z_1,z_2)\in\XX$ 
\st $z_1\neq x_1$  and $z_2\ne x_2$, the set 
\begin{equation}\label{Weylchamberdef}
\Ch_{x,z}:=\{(\sigma_{x_1,z_1}(t_1),\sigma_{x_2,z_2}(t_2)): 0\le t_1\le d_1(x_1,z_1),\, 0\le t_2 \le d_2(x_2,z_2)\}
\end{equation}
is called the {\hl Weyl chamber} with apex $x$ containing  $z$. Notice that if $z_1=x_1$ or $z_2=x_2$, then $\sigma_{x,z}$ is not defined, so the assignment in (\ref{Weylchamberdef}) 
is not well-defined. In this case we define $\Ch_{x,z}$ as follows:
$$ \Ch_{x,z}:=\bigg\{\begin{array}{lcl} \  \{(y_1, \sigma_{x_2,z_2}(t))\in\XX :0\le t\le d_2(x_2,z_2),\, y_1\in \XX_1\} & \mbox{if} & x_1=z_1\\
\ \{( \sigma_{x_1,z_1}(t),y_2)\in\XX :0\le t\le d_1(x_1,z_1),\, y_2\in \XX_2\} & \mbox{if} & x_2=z_2
\,.\end{array}$$
Similarly, for $x=(x_1,x_2)\in\XX$ and $\tilde\xi=(x_1,x_2,\theta)\in\regrand$ we call 
$$\Ch_{x,\tilde\xi}:=\{(\sigma_{x_1,\xi_1}(t_1),\sigma_{x_2,\xi_2}(t_2)): t_1, t_2\ge 0\}$$
the Weyl chamber with apex $x$ in the class of $\tilde\xi$. If $\tilde\xi\in\singrand$, we set
\begin{equation}\label{Weylsingular}
  \Ch_{x,\tilde\xi}:=\bigg\{\begin{array}{lcl}\ \{( \sigma_{x_1,\xi_1}(t),y_2)\in\XX :t\ge 0,\, y_2\in \XX_2\} & \mbox{if} & \tilde\xi=\xi_1\in (\rand)_0\cong\rand_1\,\\
 \  \{(y_1, \sigma_{x_2,\xi_2}(t))\in\XX :t\ge 0,\, y_1\in \XX_1\} & \mbox{if} & \tilde\xi=\xi_2\in (\rand)_{\pi/2}\cong \rand_2\,.\end{array}
\end{equation}
%and for $z=x$ we have $\Ch_{x,z}=\regrand$.
In this way we have defined $\Ch_{x,z}$ for any $x\in\XX$ and $z\in\overline\XX \setminus\{x\}$.
%By abuse of notation we denote for $x\in\XX$ and a subset $U\subseteq\Frand$ 
%$$\Ch_{x,U}:=\{\sigma_{x,\tilde\xi}(t): t\ge 0,\; \tilde\xi\in\preim(U) \}$$
%the set of Weyl chambers with apex $x$ and extremities in $U$.

The {\hl Weyl chamber shadow} of a set $B\subset \XX$ viewed from $x=(x_1,x_2)\in\XX\setminus B$ is defined by
\begin{equation}\label{shadowdef}
 \Sh(x:B):=\{z\in\ganz: p_1(z)\ne x_1,\, p_2(z)\ne x_2,\,  \Ch_{x,z}\cap B\ne\emptyset  \}\,.\end{equation}
It consists of all  Weyl chambers with apex $x$ which intersect $B$ non-trivially. Notice that in view of (\ref{Weylsingular}) we have 
\begin{eqnarray}\label{shadowsingular}
 && \Sh(x:B)\cap (\rand)_0 = \{\tilde \xi=\xi_1\in\rand_1 : \sigma_{x_1,\xi_1}(t)\in p_1(B) \ \,\mbox{for some}\ t\ge 0\}\,, \nonumber\\
 &&  \Sh(x:B)\cap (\rand)_{\pi/2} = \{\tilde \xi=\xi_2\in\rand_2 : \sigma_{x_2,\xi_2}(t)\in p_2(B) \ \,\mbox{for some}\ t\ge 0\}\,.\qquad \ 
\end{eqnarray}

We next fix a base point $\xo=(\xo_1,\xo_2)\in\XX$. 
For $x\in\XX$ and $r>0$ we denote by $B_x(r)$ the open ball of
radius $r$ centered at $x$. If $h\in\is(\XX_1)\times\is(\XX_2)$ is \st both projections $h_1\in\is(\XX_1)$ and $h_2\in\is(\XX_2)$ are axial, we denote $\widetilde{h^+}\in\rand$ its 
attractive fixed point, and
$\widetilde{h^-}\in\rand$ its repulsive fixed point. If for $i\in\{1,2\}$ $\ h_i^{\pm}\in\rand_i$ are the attractive and repulsive  
fixed points of $h_i$, then we get
$$\widetilde{h^\pm}=(h_1^\pm,h_2^\pm, \arctan\big( l(h_2)/l(h_1)\big)\in\regrand$$ by applying the estimate~(\ref{slopeangle}) to a point $x\in \Ax(h)$. 
Hence if $h^\pm:=\pi^F(\widetilde{h^\pm})$, we have $h^{\pm}=(h_1^\pm,h_2^\pm)$.

The following proposition states that all  Weyl chamber shadows of sufficiently large balls contain a given  open set. This will be crucial in the proof of the shadow lemma. 
Notice that our idea of proof also considerably simplifies the proof of the analogous statement for one factor (see \cite[Proposition~3.6]{MR1465601} and \cite[Lemma 3.5]{MR2290453}). 
\begin{prp}\label{largeshadows}
Suppose $g=(g_1,g_2)$ and $h=(h_1,h_2)\in \is(
\XX)$ are axial isometries \st $g_i$ and $h_i$ are independent rank one elements in $\is(\XX_i)$ for $i=1,2$. Then there exist open 
neighborhoods $U_1\subset\rand_1$, $U_2\subset\rand_2$ of $h_1^+, h_2^+$ respectively, a finite set $\Lambda\subset\Gamma$ and $c_0>0$ with the 
following properties: 

If $U:=U_1\sqcup U_2\sqcup \{(\xi_1,\xi_2,\theta)\in\regrand: \xi_1\in U_1, \xi_2\in U_2,\theta\in (0,\pi/2)\}\subset\rand$, then 
for any $r\ge c_0$ and all $y\in\XX\setminus B_\xo(r)$ there exists $\alpha\in\Lambda$ \st
$$ \alpha U\subseteq \mathrm{Sh}(y:B_\xo(r))\,.$$
\end{prp}
\prf\ For $i=1,2$ and $\eta=(\eta_1,\eta_2)\in\{g^-,g^{+},h^-,h^{+}\}$ let $U_i(\eta)\subset\ganz_i$ be an arbitrary neighborhood of $\eta_i\in\rand_i$ with $\xo_i\notin U_i(\eta)$ \st 
all $U_i(\eta)$ are pairwise disjoint in $\ganz_i$. Upon taking smaller neighborhoods, Lemma~\ref{joinrankone} provides a constant  $c>0$  \st for $i\in\{1,2\}$ any pair of points  in 
distinct neighborhoods  can be joined by a rank one geodesic $\sigma_i\subset\XX_i$ with $d(\xo_i,\sigma_i)\le c$. Moreover, according to Lemma~\ref{dynrankone} (b) there exists  a constant   
$N\in\NN$ \st for all  $\gamma=(\gamma_1,\gamma_2)\in\{g,g^{-1},h,h^{-1}\}$ and $i\in\{1,2\}$
\begin{equation}\label{factorpingpong}  
\gamma_i^N\big(\ganz_i\setminus  U_i(\gamma^{-})\big)\subseteq U_i(\gamma^+)\,.
\end{equation}
Let $y\in\XX$ arbitrary. Then one of the following cases occurs:
\begin{enumerate}
\item[1.]Case: $y_1\in \ganz_1\setminus U_1(h^+)$ and $y_2\in\ganz_2\setminus U_2(h^+)$\\
Then by~(\ref{factorpingpong}) $\,h^{-N} y \in U_1(h^-)\times U_2(h^-)$.
\item[2.]Case: $y_1\in U_1(h^+)$ and $y_2\in U_2(h^+)$\\
Since $U_i(h^+)\subset \ganz_i\setminus U_i(g^-)$, $i=1,2$, we have again by ~(\ref{factorpingpong}) $\,g^Ny \in U_1(g^+)\times U_2(g^+)$. Hence we are in Case 1 for $g^Ny$, so $h^{-N} g^Ny\in U_1(h^-)\times U_2(h^-)$.
\item[3.]Case:  $y_1\in U_1(h^+)$ and $y_2\in\ganz_2\setminus \big(U_2(h^+)\cup U_2(g^-)\big)$\\
Then $g^Ny\in U_1(g^+)\times U_2(g^+)$, which yields $h^
{-N}g^Ny \in U_1(h^-)\times U_2(h^-)$.
\item[4.]Case:  $y_1\in U_1(h^+)$ and $y_2\in\ganz_2\setminus \big(U_2(h^+)\cup U_2(g^+)\big)$\\
Then $g^{-N}y\in U_1(g^-)\times U_2(g^-)$, which gives $h^{-N}g^{-N}y \in U_1(h^-)\times U_2(h^-)$.
\item[5.]Case: $y_1\in \ganz_1\setminus \big(U_1(h^+)\cup U_1(g^-)\big)$ and $y_2\in U_2(h^+)$\\
Similarly to case 3 we obtain $h^{-N} g^Ny\in U_1(h^-)\times U_2(h^-)$.
\item[6.]Case:  $y_1\in \ganz_1\setminus  \big(U_1(h^+)\cup U_1(g^+)\big)$ and $y_2\in U_2(h^+)$\\
As in case 4 we get $h^{-N}g^{-N}y\in U_1(h^-)\times U_2(h^-)$.
\end{enumerate}
So we have shown the existence of  $\alpha=(\alpha_1,\alpha_2)\in\Lambda:=\{ h^{N}, g^{N}h^N, g^{-N}h^{N}\}$ \st  $\alpha^{-1} y \in U_1(h^-)\times U_2(h^-)$. In particular, by our choice of the neighborhoods $U_i(h^\pm)$, $i=1,2$, we have  for all $z=(z_1,z_2)\in U_1(h^+)\times U_2(h^+)$ 
$$ d_i(\sigma_{\alpha_i^{-1}y_i, z_i},\xo_i)\le c\,,\qquad i\in\{1,2\}\,.$$ 
We set $L:=\max\big\{d_i(\xo_i,\lambda_i\xo_i): i\in\{1,2\}, \lambda=(\lambda_1,\lambda_2)\in\Lambda \}$. Then for $i=1,2$ 
\be
d_i(\sigma_{y_i,\alpha_i  z_i},\xo_i)&\le & d_i(\alpha_i \sigma_{\alpha_i^{-1}y_i, z_i},\alpha_i\xo_i)+d_i(\alpha_i\xo_i,\xo_i)\\
&<&  d_i(\sigma_{\alpha_i^{-1}y_i, z_i},\xo_i) + L \le c+L\,,
\ee
which implies  $\Ch_{y,\alpha z}\cap B_\xo(\sqrt2 (c+L))\ne\emptyset$. Hence the claim follows taking $U_1:=U_1(h^+)\cap\rand_1$, $U_2:=U_2(h^+)\cap\rand_2$ and $c_0:=\sqrt2(c+L)$.\qed\\[-1mm]

Recall the definition of the Busemann function from~(\ref{buseman}). The following easy lemma relates the Busemann function of the product to the Busemann functions on the factors. 
We include a proof for the convenience of the reader.
\begin{lem}\label{RelProdFactorr}  Let $x=(x_1,x_2)$, $y=(y_1,y_2)\in\XX$, $\tilde\xi=(\xi_1,\xi_2,\theta)\in\regrand$. Then
\begin{equation}\label{busproduct}
\bs_{\tilde\xi}(x,y)=\cos\theta \cdot \bs_{\xi_1}(x_1,y_1)+\sin\theta\cdot  \bs_{\xi_2}(x_2,y_2)\,.
\end{equation}
\end{lem}
\prf\ Notice that from the definition of the Busemann functions in $\XX_1$, $\XX_2$ we have
\be \bs_{\xi_1}(x_1,y_1)&=&\lim_{s\to\infty} \big( s\cos \theta-d_1(y_1,\sigma_{x_1,\xi_1}(s\cos\theta)\big)\,,\\
\bs_{\xi_2}(x_2,y_2)&=&\lim_{s\to\infty} \big( s\sin \theta-d_2(y_2,\sigma_{x_2,\xi_2}(s\sin\theta)\big)\,.
\ee
Now 
\be && s-d(y,\sigma_{x,\tilde\xi}(s))=\frac{s^2-d(y,\sigma_{x,\tilde\xi}(s))^2}{s+d(y,\sigma_{x,\tilde\xi}(s))}\\
&&\quad=\frac{s^2\cos^2\theta-d_1(y_1,\sigma_{x_1,\xi_1}(s\cos\theta))^2}{s+d(y,\sigma_{x,\tilde\xi}(s))}+\frac{s^2\sin^2\theta-d_2(y_2,\sigma_{x_2,\xi_2}(s\sin\theta))}{s+d(y,\sigma_{x,\tilde\xi}(s))}\,,
\ee
 hence the assertion is proved if
 \be &&  \lim_{s\to\infty}\frac{s\cos\theta+d_1(y_1,\sigma_{x_1,\xi_1}(s\cos\theta))}{s+d(y,\sigma_{x,\tilde\xi}(s))}=\cos\theta\quad\an\\
 &&  \lim_{s\to\infty}\frac{s\sin\theta+d_2((y_2,\sigma_{x_2,\xi_2}(s\sin\theta))}{s+d(y,\sigma_{x,\tilde\xi}(s))}=\sin\theta\,.
 \ee
 This claim follows immediately from the triangle inequalities
\be 
&&s\cdot\cos\theta-d_1(y_1,x_1)\le d_1(y_1,\sigma_{x_1,\xi_1}(s\cdot \cos\theta))\le s\cdot\cos\theta +d_1(y_1,x_1)\,,\\
&&s\cdot\sin\theta-d_2(y_2,x_2)\le d_2(y_2,\sigma_{x_2,\xi_2}(s\cdot \sin\theta))\le s\cdot\sin\theta+d_2(y_2,x_2)\,,\\
&& s-d(y,x)\le d(y,\sigma_{x,\tilde\xi}(s)) \le s+d(y,x)\,.\hspace{6.2cm}\Box
\ee\\[1mm]
Recall that  $\rand_\theta\subset\rand$ denotes the set of points of slope $\theta\in [0,\pi/2]$. Using similar arguments as in the proof above we get the following relation for singular boundary points.
\begin{lem}\label{RelProdFactors}
 Let $x=(x_1,x_2)$, $y=(y_1,y_2)\in\XX$. Then \\
 $$ \bs_{\tilde\xi}(x,y)=\bigg\{\begin{array}{lcl}\ \bs_{\xi_1}(x_1,y_1) & \mbox{if} & \tilde\xi=\xi_1\in (\rand)_0\cong\rand_1\,\\
\ \bs_{\xi_2}(x_2,y_2) & \mbox{if} & \tilde\xi=\xi_2\in (\rand)_{\pi/2}\cong \rand_2\,.\end{array}$$
\end{lem}
To simplify notation in the sequel we further define for $b=(b_1,b_2)\in\RR^2$, $x=(x_1,x_2)$, $y=(y_1,y_2)\in\XX$, $\tilde\xi\in\rand$  the 
{\hl $b$-Busemann function} $\bs^b_{\tilde\xi} (x,y)\in\RR$ via
\begin{equation}\label{buscomb}
\bs^b_{\tilde\xi}(x,y):=\Bigg\{\begin{array}{lclcl}\   b_1\bs_{\xi_1}(x_1,y_1)&+& b_2\bs_{\xi_2}(x_2,y_2) & \mbox{if  } & \tilde\xi=(\xi_1,\xi_2,\theta)\in \regrand\,\\
   \ b_1\bs_{\xi_1}(x_1,y_1) &&  & \mbox{if } & \tilde\xi=\xi_1\in (\rand_0)\cong\rand_1\,\\
 \  && b_2\bs_{\xi_2}(x_2,y_2) & \mbox{if } & \tilde\xi=\xi_2\in (\rand)_{\pi/2}\cong \rand_2\,.\end{array}
\end{equation}
%Then for $\tilde\xi=(\xi_1,\xi_2,\theta)\in\regrand$ equation~(\ref{busproduct}) translates to 
%$$ \bs_{\tilde\xi}(x,y)=\langle \bv_{\pi^F(\tilde\xi)}(x,y),H_\theta\rangle\,.$$
For convenience we denote
\begin{equation}\label{thetavector}
H_\theta:= \left(\begin{array}{c}\cos\theta\\\sin\theta\end{array}\right)\in\RR^2
\end{equation}
the unique unit vector of direction $\theta\in [0,\pi/2]$, and $\langle\cdot,\cdot\rangle$ the Euclidean inner product in $\RR^2$.  
In the sequel we will need the following 
\begin{df}\label{dirdist} The {\hl directional distance} of the ordered pair $(x, y)\in
\XX\times \XX$ \wrt the slope $\theta$ is defined by
\begin{eqnarray*}
 \bs_{\theta}:\quad\XX\times \XX &\to & \RR\\
 (x, y) &\mapsto &  \bs_{\theta} (x, y)\;:=\, \langle H_\theta,H(x,y)\rangle \,.
% \max \{ \bs_{\tilde \xi} (x,y): \tilde\xi\in\rand_\theta\}\,.
\end{eqnarray*} 
\end{df}
In particular, if $\theta=0$, then $\bs_\theta(x,y)=d_1(p_1(x),p_1(y))$,  if $\theta=\pi/2$, then $\bs_\theta(x,y)=d_2(p_2(x),p_2(y))$. 

By $\big(\is(\XX_1)\times\is(\XX_2)\big)$-invariance of the distance vector we immediately get that 
$$ \bs_\theta(gx,gy)=\bs_\theta(x,y)$$ for any $x,y\in\XX$ and $g\in  \is(\XX_1)\times\is(\XX_2)$. Moreover, the symmetry and triangle inequality for the distances $d_1$ and $d_2$ 
directly imply the symmetry and triangle inequality for $\bs_\theta$.
The following important proposition states that for $\theta\in (0,\pi/2)$ the directional distance $\bs_\theta$ is in fact a distance.
\begin{prp}\label{dirisdistance}
For  $\theta\in (0,\pi/2)$ the directional distance $\bs_\theta$ is a distance.
\end{prp}
\prf \ Let $x=(x_1,x_2)$, $y=(y_1,y_2)\in\XX$. We clearly have 
$$\bs_\theta(x,y)=\cos\theta\cdot d_1(x_1,y_1)+\sin\theta\cdot d_2(x_2,y_2)\ge 0\,,$$ because all terms involved are non-negative. 
Moreover, if $\bs_\theta(x,y)=0$, then $\cos\theta>0$ and $\sin\theta>0$ imply $ d_1(x_1,y_1)=0$ and $ d_2(x_2,y_2)=0$, hence  $x=y$.

Finally, we have already noticed that the symmetry and triangle inequality follow directly from the symmetry and triangle inequality for the distances $d_1$ and $d_2$. \qed\\[-1mm]

The following easy facts will be convenient in the sequel.
\begin{lem} \label{dirBus}
Let $x\in\XX$  and $\tilde\xi\in \rand_\theta$ for some $\theta\in [0,\pi/2]$. Then 
$$ y\in \Ch_{x,\tilde\xi}\setminus \{x\}\ \iff\quad 
\bs_\theta(x,y)=\bs_{\tilde\xi}(x,y)\,.$$
\end{lem}
\prf Write $x=(x_1,x_2)$ and let  $y=(y_1,y_2)\in\XX\setminus \{x\}$  arbitrary. Lemma~\ref{RelProdFactorr}, Lemma~\ref{RelProdFactors} and (\ref{bounded}) imply 
\begin{equation}\label{ineqbusedir}
\bs_{\tilde\xi}(x,y)\le\cos\theta \cdot d_1(x_1,y_1)+\sin\theta\cdot  d_2(x_2,y_2)=\langle H_\theta, H(x,y)\rangle\,.
\end{equation} 
Assume first that $\theta\in (0,\pi/2)$ and write $\tilde \xi=(\xi_1,\xi_2,\theta)$. Then we have equality in (\ref{ineqbusedir}) if and only if $\bs_{\xi_i}(x_i,y_i)=d_i(x_i,y_i)$ for $i=1,2$.  If $x_1\neq y_1$ and  $x_2\neq y_2$, this is precisely the case if $y_i$ is a point on the geodesic ray $\sigma_{x_i\xi_i}$ for $i=1,2$ which is equivalent to $y\in  \Ch_{x,\tilde\xi}$. If $x_1=y_1$, then 
$\bs_{\tilde\xi}(x,y)=\cos\theta\cdot  0 +\sin\theta\cdot d_2(x_2,y_2)$ if and only if $y_2$ is a point on the geodesic ray $\sigma_{x_2,\xi_2}$. This again holds if and only if $y\in\Ch_{x,\tilde\xi}$. The case $x_2=y_2$ is analogous.

If $\tilde\xi=\xi_1\in (\rand)_{0}\cong\rand_1$  we have  equality in (\ref{ineqbusedir}) if and only if $y_1$ is a point on the geodesic ray $\sigma_{x_1,\xi_1}$. This is equivalent to $y\in\Ch_{x,\tilde\xi}$. 

Similarly, if  $\tilde\xi=\xi_2\in(\rand)_{\pi/2}\cong\rand_2$ we have  equality in (\ref{ineqbusedir}) if and only if $y\in\Ch_{x,\tilde\xi}$. \qed\\[-1mm]

If $\XX$ is geodesically complete, this lemma allows to give the following nice geometric interpretation of the directional distance.
%, which will be convenient in the proof of Proposition~\ref{dirisdistance}.
\begin{cor}\label{dirdistbusemann}
If $\XX=\XX_1\times \XX_2$ is geodesically complete, $\theta\in [0,\pi/2]$, $x=(x_1,x_2)$, $y=(y_1,y_2)\in\XX$, then %If for $i\in\{1,2\}$  
%$x_i=y_i$ or the geodesic $\sigma_{x_i,y_i}\subset \XX_i$ can be extended to a geodesic ray, then
$$\bs_\theta(x,y)= \max \{ \bs_{\tilde \xi} (x,y): \tilde\xi\in\rand_\theta\}\,.$$
If $\theta=0$, then the conclusion holds under the weaker condition that $\XX_1$ is geodesically complete, 
%$x_2=y_2$ or the geodesic  $\sigma_{x_2,y_2}\subset \XX_2$ can be extended to a geodesic ray, 
if $\theta=\pi/2$, the conclusion holds under the condition that $\XX_2$ is geodesically complete. % $x_1=y_1$ or the geodesic  $\sigma_{x_1,y_1}\subset \XX_1$ can be extended to a geodesic ray. 
\end{cor}
\prf If $y=x$, then $\bs_\theta(x,y)=\langle H_\theta,H(x,y)\rangle=0$ and 
$\bs_{\tilde\xi}(x,y)\le d(x,y)=0$ for all $\tilde\xi\in\rand_\theta$. Hence we have $\max \{ \bs_{\tilde \xi} (x,y): \tilde\xi\in\rand_\theta\}=0=\bs_\theta(x,y)$. 

We next treat the case $y\ne x$. 
Since $\XX$ is geodesically complete, every point $y\in\XX$ is contained in a Weyl chamber $\Ch_{x,\tilde\xi}$ for some $\tilde\xi\in\rand_\theta$. Hence if $y\ne x$, the previous lemma 
implies $\bs_\theta(x,y)=\bs_{\tilde\xi}(x,y)$. Moreover, if $\tilde\zeta\in\rand_\theta$ is arbitrary, then by Lemma~\ref{RelProdFactorr}, Lemma~\ref{RelProdFactors} and (\ref{bounded})
$$\bs_{\tilde\zeta}(x,y)\le \langle H_\theta,H(x,y)\rangle\,.$$ 
Summarizing we conclude $\quad\bs_\theta(x,y)=\bs_{\tilde\xi}(x,y)=\max \{ \bs_{\tilde \xi} (x,y): \tilde\xi\in\rand_\theta\}\,.$

In order to prove the remaining assertions we write   $x=(x_1,x_2)$, $y=(y_1,y_2)$ and assume that $y\ne x$.
First assume that $\theta=0$ and $\XX_1$ is geodesically complete. Then there exists $\xi_1\in\rand_1$ \st $y_1=\sigma_{x_1,\xi_1}(t)$ for some $t\ge 0$. Hence if 
$\tilde\xi\in(\rand)_0$ 
is the unique point identified with $\xi_1\in\rand_1$ we have $y\in\Ch_{x,\tilde\xi}\setminus\{x\}$. The claim now follows as before from the previous lemma, Lemma~\ref{RelProdFactors} and 
(\ref{bounded}). 

The case $\theta=\pi/2$ and $\XX_2$ geodesically complete is analogous.\qed\\[-1mm]

Recall the definition of Weyl chamber shadows from~(\ref{shadowdef}). The following lemma will be needed in the proof of the shadow lemma Theorem~\ref{shadowlemma}.
\begin{lem}\label{esti}
Let $c>0$, $z=(z_1,z_2) \in\XX$ with $d(o,z)>c$, and $\tilde\eta\in\Sh\big(o:B_{z}(c)\big)\cap\rand$. If $\tilde\eta=(\eta_1,\eta_2,\theta)\in\regrand$, then 
$$  0\le d_1(o_1,z_1)-\bs_{\eta_1}(o_1,z_1)<2c\quad\an\\ \quad 0\le d_2(o_2,z_2)-\bs_{\eta_2}(o_2,z_2)<2c\,,$$
if $\tilde\eta=\eta_1\in (\rand)_0\cong \rand_1$, then  $\quad 0\le d_1(o_1,z_1)-\bs_{\eta_1}(o_1,z_1)<2c$, and \\
if $\tilde\eta=\eta_2\in (\rand)_{\pi/2}\cong \rand_2$, then  $\quad 0\le d_2(o_2,z_2)-\bs_{\eta_2}(o_2,z_2)<2c$.
\end{lem} 
\prf\  By definition $\tilde\eta\in \Sh\big(o:B_{z}(c)\big)$ translates to $\Ch_{\xo,\tilde\eta}\cap B_z(c)\ne\emptyset$. Hence if $\tilde\eta=(\eta_1,\eta_2,\theta)\in\regrand$  there exist 
$t_1,t_2\ge 0$  \st $d\big((z_1,z_2), (\sigma_{\xo_1,\eta_1}(t_1),\sigma_{\xo_2,\eta_2}(t_2))\big)<c$. Therefore we have $d_i(z_i,\sigma_{\xo_i,\eta_i}(t_i))<c $ for $i\in\{1,2\}$, so the 
claim follows from Lemma~\ref{distbusbounded}. The conclusion for $\tilde \eta\in\singrand$ is clear in view of Lemma~\ref{RelProdFactors}, (\ref{shadowsingular}) and Lemma~\ref{distbusbounded}. \qed

\section{The exponent of growth}\label{ExpGrowth}

For the remainder of the article  $\XX$ is a product of locally compact 
Hadamard spaces $\XX_1$, $\XX_2$, $\xo=(\xo_1, \xo_2)$ a fixed base point, and $\Gamma\subset\is(\XX_1)\times\is(\XX_2)$ a discrete group which contains two isometries $g=(g_1,g_2)$ and 
$h=(h_1,h_2)$ \st for $i=1,2$ $g_i$ and $h_i$ are independent rank one elements of $\Gamma_i$. 
Recall that the geometric limit set of a  group $\Gamma$ acting  by isometries on a locally compact Hadamard space is defined by
$\Lim:=\overline{\Gamma\at x}\cap\rand$, where $x\in \XX$ is arbitrary.
In this section we recall the notion of  exponent of  growth introduced in \cite{MR2629900} and give an important criterion for divergence or convergence of certain sums over $\Gamma$. 
This will play a central role in the construction of (generalized) Patterson-Sullivan measures in Sections~\ref{ClasPatSulconst} and \ref{genPseries}.

We recall the notation introduced in Section~\ref{prodHadspaces} and put for $x,y\in\XX$, $\theta\in[0,\pi/2]$, $\eps>0$ 
$$ \Gamma(x,y;\theta,\eps):=\{\gamma\in\Gamma: \gamma y\ne x\quad \mbox{and}\ \ |\theta(x,\gamma y)-\theta | <\eps \}\,.$$
In order to define the exponent of growth of $\Gamma$ of slope $\theta$ we put
$$ \delta_\theta^\eps(x,y):=\inf \{s>0: \sum_{\gamma\in
  \Gamma(x,y;\theta,\eps)}
  e^{-s d(x,\gamma y)} \ \ \mbox{converges}\}\,.$$
If $\delta(\Gamma)$ denotes the {\hd critical exponent} of $\Gamma$ defined by
\begin{equation}\label{critexp}
  \delta(\Gamma):=\inf\{s>0: \sum_{\gamma\in\Gamma}e^{-sd(\xo,\gamma\xo)}\ \mbox{converges}\}\,,
\end{equation}
we clearly have $\delta_\theta^\eps(x,y)\le \delta(\Gamma)$ with equality if $\eps>\pi/2$.
Moreover, Lemma~6.1 in \cite{MR2629900} shows that  $\delta_\theta^\eps(x,y)$ is related to the numbers  
$$\Delta N_\theta^\eps(x,y;n):=\#\{ \gamma\in\Gamma\;:\, n-1< d(x,\gamma y)\le n\,,\ |\theta(x,\gamma y)-\theta|<\eps\}\,,\quad n\in\NN \,,$$
via 
\begin{equation}\label{defbylimsup}
\delta_{\theta}^{\eps}(x,y)=\limsup_{n\to\infty}\frac{\log
  \Delta N_{\theta}^{\eps}(x,y;n)}{n}
  \end{equation}
and thus can be interpreted as %the fastest 
an  exponential growth rate of the number of orbit points with slope $\eps$-close to $\theta$.

Recall that the {\hl exponent of growth of $\Gamma$ of slope $\theta$} is  defined by 
$$\delta_{\theta}(\Gamma):=\liminf_{\eps\to
  0}\delta_{\theta}^{\eps}(o,o)\,.$$ 
 Notice that this number  $\delta_\theta(\Gamma)$  does not depend on the choice of arguments 
of $\delta_{\theta}^\eps$ by Lemma~6.3 in \cite{MR2629900}, and 
$\delta_\theta(\Gamma)\le \delta(\Gamma)$ for all $\theta\in [0,\pi/2]$. 
%Moreover, in the definition of $\delta_\theta(\Gamma)$ for $\theta\in (0,\pi/2)$ one may substitute $$\#\{ \gamma\in\Gamma\;:\, \gamma y\ne x\,,\ d(x,\gamma y)<R\,,\ \Big| \frac{d_2(p_2(\gamma y),p_2(x))}{d_1(p_1(\gamma y),p_1(x))}-\tan\theta\Big|<\eps\}\,$$ in~(\ref{defbylimsup}) instead of $N_\theta^\eps(x,y;R)$. 

Furthermore, we recall the following properties from Section 6 in \cite{MR2629900}:\\[3mm]
\begin{pros}
\begin{enumerate}
\item[(a)] If $\Lim\cap  \rand_\theta\ne \emptyset$, then $\delta_{\theta}(\Gamma)\ge 0$.
\item[(b)] If  $(\theta_j)\subset [0,\pi/2]$ is a sequence converging to $\theta\in [0,\pi/2]$, then
$$ \limsup_{j\to\infty}\,\delta_{\theta_j}(\Gamma)\le \delta_{\theta}(\Gamma)\,.$$
\end{enumerate}
\end{pros}

It will turn out useful to  extend the exponent of growth to a homogeneous map $\Psi_\Gamma:\RR_{\ge 0}^2\to \RR$ as follows: If $x=(x_1,x_2)\in\RR_{\ge 0}^2$ we put $\theta(x):=\arctan (x_2/x_1)$ 
and set
\begin{equation}\label{defPsi}
\Psi_\Gamma(x):=||x||\cdot \delta_{\theta(x)}\,.
\end{equation}
In \cite{MR2629900} we showed that $\Psi_\Gamma$ is concave. This implies in particular that there exists a unique $\theta_*\in [0,\pi/2]$ such that 
$\delta_{\theta_*}(\Gamma)=\max\{\delta_\theta(\Gamma):\theta\in [0,\pi/2]\}$. The following important proposition will play a key role in the proof of Theorem~A and for the construction of generalized 
Patterson-Sullivan measures.
%This property will be essential for the proof of Theorem~A and Proposition~\ref{key}, where we describe the construction of generalized conformal densities. \\[1mm]
Recall the definition of the distance vector  and of $H_\theta$ from~(\ref{distancevector}) and~(\ref{thetavector}) respectively.
\begin{prp}\label{convdiv}
Let $D\subseteq [0,\pi/2]$ be a relatively  open interval, set 
$\Gamma_D:=\{\gamma\in\Gamma : \theta(\xo,\gamma\xo)\in D\}$, and let $f: \RR_{\ge 0}^2\to \RR $ be a continuous homogeneous function.
\begin{enumerate}
\item[(a)]
If there exists $\hat\theta\in D$ \st $f(H_{\hat\theta})<\delta_{\hat\theta}(\Gamma)\,,$
then the series $\sum_{\gamma\in\Gamma_D}e^{-f(H(\xo,\gamma\xo))}$ diverges.
\item[(b)]
If $f(H_\theta)> \delta_{\theta}(\Gamma)$ 
for all $\theta\in \overline{D}$,
 then the series  $\sum_{\gamma\in\Gamma_D}e^{-f(H(\xo,\gamma\xo))}$ converges.
\end{enumerate} 
\end{prp}
\prf\  For $\gamma\in\Gamma$ we set $H_\gamma:=H(\xo,\gamma\xo)/d(\xo,\gamma\xo)$.

(a) Let $\hat\theta \in D$  \st 
$\  f(H_{\hat\theta})<\delta_{\hat\theta}(\Gamma)\,.$  
Since $\delta_{\hat\theta}(\Gamma)=\liminf_{\eps\to
  0}\delta_{\hat\theta}^\eps(o,o)$, there exists
$\eps\in (0,\pi/4)$ and $ \hat s\in\RR$ \st for 
$\gamma\in\Gamma_D$ with  
$\mid \theta(\xo, \gamma \xo)-\hat\theta \mid <\eps$  we have  
$$  f(H_\gamma)< \hat s<\delta_{\hat\theta}^\eps(o,o)\,.$$
Therefore  
$$ \sum_{\gamma\in\Gamma_D}e^{-f(H(\xo,\gamma\xo))}
> \sum_{\begin{smallmatrix}{\gamma\in\Gamma}\\{|\theta(\xo,\gamma o)-\hat\theta|<\eps}\end{smallmatrix}}e^{- \hat s d(o,\gamma o)},$$
and the latter sum diverges since $ \hat s<\delta_{\hat\theta}^\eps(o,o)$.
\smallskip

(b)
Let $ \hat\theta\in\overline{D}$. Since $f(H_{\hat\theta})> \delta_{\hat\theta}(\Gamma)=\liminf_{ \eps\to 0}\delta_{\hat\theta}^ \eps(\xo,\xo)$, there exists  $ \eps'\in (0,\pi/4)$ and
  $ \hat s<f(H_{\hat\theta})$ \st 
\begin{equation}\label{sandwich} 
\delta_{  \hat\theta}^{ \eps'}(o,o)< \hat s< f(H_{\hat\theta})\,.
\end{equation}
For $\theta\in [0,\pi/2]$ and $\eps>0$ we put $B_{\theta}( \eps):=\{\theta'\in [0,\pi/2] : |\theta'-\theta|< \eps\}$. The continuity of the function $f$
and inequality~(\ref{sandwich}) imply
the existence of $ \hat\eps< \eps'$ \st for any $\theta\in
B_{  \hat\theta}( \hat\eps)$  we have $ \hat s<f(H_\theta)$. Hence for all $z\in\XX$ with $\theta(o,z)\in B_{  \hat\theta}( \hat\eps)$ we have
$$\frac{f(H(\xo,z))}{d(\xo,z)}>\hat s> \delta_{  \hat\theta}^{ \eps'}(o,o) \ge \delta_{  \hat\theta}^{ \hat \eps}(o,o)\,.$$
We now choose a sequence $(\theta_j)\subset\overline{D}$ and corresponding
sequences $(s_j)\subset \RR^+$ and $( \eps_j)\subset\RR^+$ \st for every 
$\theta\in B_{ \theta_j}( \eps_j)$   we have 
$$\delta_{ \theta_j}^{ \eps_j}(o,o)<s_j<f(H_\theta)\,,\ \an\quad \overline{D}\subseteq \bigcup_{j\in\NN} B_{\theta_j}(\eps_j)\,.$$ 
Since $\overline{D}$ is compact we may
extract a finite covering $\bigcup_{j=1}^l B_{ \theta_j}( \eps_j)\,$,
and conclude 
\begin{align*}
\sum_{\gamma\in\Gamma_D}e^{-f(H(\xo,\gamma\xo))}& \le
\sum_{j=1}^l 
\sum_{\begin{smallmatrix}{\gamma\in\Gamma}\\{|\theta(\xo,\gamma o)-
     \theta_j|< \eps_j}\end{smallmatrix}} e^{-f(H(\xo,\gamma\xo))}\\
&\le\sum_{j=1}^l \sum_{\begin{smallmatrix}{\gamma\in\Gamma}\\{|\theta(\xo,\gamma o)-
     \theta_j|< \eps_j}\end{smallmatrix}}e^{-s_j d(o,\gamma o)} <\infty\,,
\end{align*}
because $s_j>\delta_{ \theta_j}^{ \eps_j}(o,o)$ for $1\le j\le l$.\qed\\[-1mm]

Taking $D=[0,\pi/2]$ and $f(H)=s\cdot \Vert H\Vert$ we obtain as a corollary that $\delta(\Gamma)=\max\{\delta_\theta(\Gamma):\theta\in [0,\pi/2]\}=\delta_{\theta_*}(\Gamma)$. 
 We conclude this section with two illustrative examples. \\[0mm]

\noindent{\sc Example 1} (see \cite[Section 6]{MR2629900}) $\quad$
If $\XX$ is a product $\XX=\XX_1\times \XX_2$ of Hadamard manifolds with pinched negative curvature, and 
$\Gamma_1\subset\is(\XX_1)$, $\Gamma_2\subset\is(\XX_2)$  are convex cocompact
groups with critical exponents $\delta_1, \delta_2$, then for $\Gamma:=\Gamma_1\times\Gamma_2$ and for every $\theta\in [0,\pi/2]$
$$ \delta_\theta(\Gamma)=\delta_1\cos\theta +\delta_2\sin\theta\,.$$
This number is maximal for $\theta_*=\arctan(\delta_2/\delta_1)$ and we have $\delta(\Gamma)=\delta_{\theta_*}(\Gamma)=\sqrt{\delta_1^2+\delta_2^2}$. 
The homogeneous function $\Psi_\Gamma:\RR_{\ge 0}^2\to \RR$ is simply the linear functional defined by $\Psi_\Gamma=\langle \left(\ba{c} \delta_1\\\delta_2\ea\right),\cdot\rangle$. \\[2mm]

\noindent{\sc Example 2}  $\quad$  Consider a product of hyperbolic planes $\XX=\HH^2\times \HH^2$ and a Hilbert modular group $\Gamma\subset\is(\XX)$. Then $\Gamma$ is an irreducible 
non-uniform lattice in a higher rank symmetric space,
hence from Proposition 7.2 and 7.3 in \cite{MR1675889} we know that $\Psi_\Gamma=\langle \left(\ba{c}1\\1\ea\right),\cdot\rangle$, i.e.
$\delta_\theta(\Gamma)=\cos\theta +\sin\theta$. Here $\delta_\theta(\Gamma)$ is maximal for $\theta_*=\pi/4$ and we have $\delta(\Gamma)=\delta_{\theta_*}(\Gamma)=\sqrt2$.

\section{The classical Patterson-Sullivan construction}\label{ClasPatSulconst}

In this section we will  construct  a conformal density for $\Gamma$ using an idea originally due to  S.~J.~Patterson (\cite{MR0450547}) in the context 
of Fuchsian groups. Taking advantage of Proposition~\ref{convdiv} we will be able to describe precisely its support and hence 
 prove Theorem~A. 

Recall that a $\Gamma$-invariant conformal density of dimension $\alpha\ge 0$ is  a
continuous map $\mu$ from $\XX$ to the cone $\MM^+(\rand)$ of positive   finite Borel measures
on $\rand$ \st 
$\supp(\mu_{\xo})\subseteq\Lim$, $\gamma*\mu_x=\mu_{\gamma^{-1}x}$ for all $\gamma\in\Gamma$, $x\in\XX$ and
$$ \frac{d\mu_x}{d\mu_\xo}(\tilde\eta)=e^{\alpha \bs_{\tilde \eta}(\xo,x)}\qquad\mbox{for}\ \tilde\eta\in\supp(\mu_{\xo}),\, x\in\XX\,.$$

In order to construct a $\Gamma$-invariant conformal density of dimension $\delta(\Gamma)$ we first suppose that we are given a map $b:\Gamma\to\RR$, $\gamma\mapsto b_\gamma$, \st the sum 
$$\sum_{\gamma\in\Gamma} e^{-sb_\gamma}$$ converges for $s>1$ and diverges for $s<1$. The following useful lemma states that if the above sum converges for $s=1$, then we can slightly modify it to 
obtain a sum which diverges for $s\le 1$ and converges for $s>1$.   

\begin{lem}[\ (Patterson \cite{MR0450547}\label{Pt})]
There exists a positive increasing function $h$ on $[0,\infty)$ \st 
\begin{enumerate}
\item[\rm(i)]
$\sum_{\gamma\in \Gamma} e^{-s b_\gamma}
h(b_\gamma)$ has exponent of convergence $s=1$ and
diverges at $s=1$;
\item[\rm (ii)]
for any $\varepsilon>0$ there exists $r_0>0$  \st  for $r\ge r_0$ and $t>1$
 $$h(rt)\le t^\varepsilon h(r)\,.$$
\end{enumerate}
\end{lem}

Recall the definition of the exponent of growth of $\Gamma$ and its properties from Section~\ref{ExpGrowth}. We have already noticed that there exists a unique 
$\theta_*\in [0,\pi/2]$ \st $\delta(\Gamma)=\delta_{\theta_*}(\Gamma)$. 

Following the original idea of Patterson \cite{MR0450547}, we apply the above lemma to the map 
\begin{equation}\label{bClas}
b:\Gamma\to\RR, \ \gamma\mapsto \delta(\Gamma)\cdot d(\xo,\gamma\xo)\,. 
\end{equation}
Then by definition~(\ref{critexp}) of the critical exponent $\delta(\Gamma)$ the series 
$\ \,\sum_{\gamma\in\Gamma}e^{-sb_\gamma}\ $ has exponent of convergence $s=1$. 
If this sum converges at $s=1$ we take the increasing function $h$ from the previous lemma, otherwise we set $h\equiv 1$ and define
$$\Summ^s:=\sum_{\gamma\in\Gamma}e^{-sb_\gamma}h(b_\gamma)\,.$$
We then obtain a family of orbital measures on $\overline{\XX}$ as follows:
If  $D$ denotes the unit Dirac point measure, then for $x \in\XX$ and $s>1$ we set
$$ \mu_x^s:=\frac1{\Summ^s}\sum_{\gamma\in\Gamma}e^{-s\delta(\Gamma) d(x,\gamma\xo)}h(b_\gamma)D(\gamma\xo)\,.$$
These measures are $\Gamma$-equivariant by construction and absolutely
continuous \wrt each other. 

Let $(\Cnt^{0}(\ganz),\Vert\cdot\Vert_\infty)$ denote the
space of real
valued continuous functions on $\ganz$ with norm
$\Vert f\Vert_\infty=\max\{|f(x)| :  x\in\ganz\}$,
$f\in \Cnt^{0}(\ganz)$. 
We endow  the cone $\MM^+ (\overline\XX)$ of positive finite Borel
measures on $\overline{\XX}$  
with the pseudo-metric 
\begin{equation}\label{pseudometricmeasures}
\rho (\mu_1,\mu_2):=\sup\bigg\{\bigg|\int_{\ganz}
f\;d\mu_1 - \int_{\ganz} f\;d\mu_2 \bigg| :  f\in
\Cnt^{0}(\ganz)\,,\ \Vert f\Vert_\infty =1\bigg\}\,,
\end{equation}
$\mu_1,\mu_2\in \MM^+ (\overline\XX)$, 
and obtain the
following:

\begin{lem}\label{equiconclass} 
The family of maps 
${\mathcal F}:=\{x\mapsto\mu_x^s :  1<s
\le 2\}$ from $\XX$ to $\MM^+ (\overline\XX)$ is equicontinuous. 
\end{lem}
\prf\  Let $x,y\in\XX$.  For $\gamma\in\Gamma$ we abbreviate
$$ q_\gamma(y,x):=\delta(\Gamma)\big(d(y,\gamma\xo)-d(x,\gamma\xo)\big)$$
and notice that
\begin{equation}\label{estimq}
|q_\gamma(y,x)|\le \delta(\Gamma)d(x,y)\,.
\end{equation}
If $s\in(1, 2]$ and $f\in\Cnt^{0}(\ganz)$, the inequality $|1-e^{-t}|\le
e^{|t|}-1$, $t\in\RR$, then gives 
\begin{align*}
\bigg|\int_{\ganz}
f\;d\mu_x^s - \int_{\ganz} f\;d\mu_y^s\bigg|&\le\frac1{
\Summ^s}\sum_{\gamma\in \Gamma}
e^{-s\delta(\Gamma)d(x,\gamma\xo)}\cdot h(b_\gamma) |f(\gamma\xo) |\,\big|1- e^{-s q_\gamma(y,x)}\big|\\
 &\le\frac{\Vert f\Vert_\infty}{\Summ^s}\sum_{\gamma\in 
\Gamma}e^{-s b_\gamma}h(b_\gamma) e^{-sq_\gamma(x,o)} (e^{s
|q_\gamma(y,x)|}-1)\,.\end{align*}
Since $f\in\Cnt^{0}(\ganz)$ was arbitrary, $s\le 2$ and
$\sum_{\gamma\in\Gamma} e^{-sb_\gamma}h(b_\gamma)=\Summ^s$,
we conclude using 
(\ref{estimq}) 
$$
\rho(\mu_x^s,\mu_y^s)\le e^{2\delta(\Gamma)d(o,x)}\at\big
(e^{2\delta(\Gamma)d(x,y)}-1\big)\,.$$
This proves that ${\mathcal F}$ is equicontinuous.\qed
\begin{lem}
For any $x\in\XX$ there exists a sequence 
\hbox{$(s_n)\searrow 1$} \st the measures
 $\mu_{x}^{s_n}\subset \MM^+(\overline\XX)$ converge weakly to a measure $\mu_{x}:=\mu_x( \theta,\tau,b)$ as $n\to\infty$. \end{lem}
\prf\   The compactness of the space $\ganz$ implies that every sequence of measures in $\MM^+(\overline\XX)$ possesses
a weakly convergent subsequence. \qed\\[-1mm]

The Theorem of Arzel{\`a}-Ascoli \cite[Theorem 7.17, 
p.$\,$233]{Kelley} now allows  to conclude that
${\mathcal F}$ is relatively compact in the space of continuous maps $\Cnt(\XX,\MM^+
(\overline\XX))$ endowed with the topology of uniform convergence on compact
sets. 
From the definition of $(\mu_{x}^{s})_{x\in\XX}$ it
follows that every accumulation point
$\mu=(\mu_x)_{x\in\XX}$ of 
${\mathcal F}$ 
as $s\searrow 1$  takes its values in  $\MM^+(\rand)$. Moreover, the following proposition shows that the family of measures obtained in this way are absolutely continuous with Radon-Nikodym derivative
$$\frac{d\mu_x}{d\mu_\xo}(\tilde \eta)=e^{\delta(\Gamma)\bs_{\tilde\eta}(\xo,x)}$$
for $\tilde\eta\in \supp(\mu_\xo)$, $x\in\XX$.
\begin{prp} Every accumulation point $\mu=(\mu_x)_{x\in\XX}$ of the
family ${\mathcal F}$ in\\
$\Cnt(\XX,\MM^+(\overline\XX))$ is a $\delta(\Gamma)$-dimensional conformal density. 
\end{prp}
\prf\ Let $(\mu_x)_{x\in\XX}$ be an accumulation point of ${\mathcal
 F}$. By construction, the measures $\mu_x$, $x\in\XX$, are
 $\Gamma$-equivariant and supported on the limit set $\Lim$. 
It therefore remains to show  
$$\frac{d\mu_x}{d\mu_\xo}(\tilde\eta)=e^{\delta(\Gamma)\bs_{\tilde\eta}(\xo,x)}\quad\mbox{for any}\ x\in\XX\,,\ \tilde\eta\in\supp(\mu_\xo)\,.$$
Notice that if $(y_n)\subset \XX$ is a sequence converging to 
a point $\tilde\eta\in \rand$,  then $d(x,y_{n})-d(\cdot, y_{n})\to \bs_{\tilde\eta}(x,\cdot)$ uniformly on compact sets in $\XX$. Hence it suffices to prove that
for any $f\in\Cnt^0(\ganz)$
$$ \lim_{s\searrow 1}\bigg|\int_{\ganz}f(z)d\mu_\xo^s(z)-\int_{\ganz}f(z)e^{-\delta(\Gamma)( d(\xo,z)-d(x,z))} d\mu_x^s(z)\bigg|=0\,.$$
We choose $f\in \Cnt^0(\ganz)$ arbitrary, $s\in (1,2]$ and estimate  
\begin{eqnarray*}
\bigg|\int_{\ganz} f(z)d\mu_\xo^s(z)& - & \int_{\ganz}f(z)e^{-\delta(\Gamma)( d(\xo,\gamma\xo)-d(x,\gamma\xo))}d\mu_x^s(z) \bigg|\\
&=& \frac1{\Summ^s}\bigg| \sum_{\gamma\in\Gamma} f(\gamma\xo)h(b_\gamma)\left( e^{-sb_\gamma}-  e^{-\delta(\Gamma)( d(\xo,\gamma\xo)-d(x,\gamma\xo))}e^{-s\delta(\Gamma)d(x,\gamma\xo)}\right)\bigg|\\
& \le & \frac1{\Summ^s}\sum_{\gamma\in\Gamma}|f(\gamma\xo)|
\cdot e^{-sb_\gamma}h(b_\gamma)\cdot \big|1-e^{\delta(\Gamma)(s-1)( d(\xo,\gamma\xo)-d(x,\gamma\xo))}\big|\\
&\le &  \frac{\Vert f\Vert_\infty}{\Summ^s} \sum_{\gamma\in\Gamma}  e ^{-sb_\gamma}h(b_\gamma)\cdot\, \big(e^{\delta(\Gamma)(s-1) |d(\xo,\gamma\xo)-d(x,\gamma\xo)|}-1\big)\\
&\le & \Vert f\Vert_\infty\cdot \big( e^{\delta(\Gamma)(s-1)d(\xo,x)}-1\big)\,. 
\end{eqnarray*}
Since the last term tends to zero as $s$ tends to $1$ the claim follows.\qed\\[-1mm]

Recall that $\theta_*\in [0,\pi/2]$ is the unique point \st $\delta_{\theta_*}(\Gamma)=\delta(\Gamma)$. 
In order to prove Theorem A it remains to show that the support of the conformal density $\mu$ constructed  above is included in the unique 
$\Gamma$-invariant subset of the limit set which consists of all limit points with slope 
$\theta_*$. This is the content of the following
\begin{prp}
The support of the conformal density $\mu=(\mu_x)_{x\in\XX}$ is contained in $\Lim\cap\rand_{\theta_*}$. 
\end{prp}
\prf\ Recall the definition of $b_\gamma$ from~(\ref{bClas}) and let $h$ be the function from Lemma \ref{Pt}. Using Proposition~\ref{convdiv} we will prove that for any $\eps>0$ 
$$\sum_{\begin{smallmatrix}{\gamma\in\Gamma}\\
{|\theta(\xo,\gamma\xo)-\theta_*|>\eps}\end{smallmatrix}} e^{-b_\gamma} h(b_\gamma) <\infty\,.$$
Let $\eps >0$ arbitrary and set 
$$s_\eps:=\max\{\delta_\theta(\Gamma): \theta\in [0,\pi/2], \, |\theta-\theta_*|\ge \eps\}\,.$$
Then by choice of $\theta_*$ we have $\delta(\Gamma)=\delta_{\theta_*}(\Gamma)>s_\eps$. Fix $\alpha:=\frac12\big(\delta(\Gamma)-s_\eps\big)$ and let $r_0>0$ \st for all $r\ge r_0$ 
and $t>1$ we have $h(rt)\le t^\alpha h(r)$. In particular, if 
$b_\gamma\ge r_0$, then 
$$h(b_\gamma)=h(\frac{b_\gamma}{r_0}\cdot r_0)\le \left(\frac{b_\gamma}{r_0}\right)^\alpha\cdot h(r_0)=\frac{\delta(\Gamma)^\alpha h(r_0)}{r_0^\alpha}\cdot e^{\alpha\log d(\xo,\gamma\xo)}\,.$$
Set $\Gamma_\eps:=\{\gamma\in\Gamma: |\theta(\xo,\gamma\xo)-\theta_*|>\eps,\, b_\gamma\ge r_0\}\,.$ Then 
\be \sum_{\gamma\in\Gamma_\eps}e^{-b_\gamma} h(b_\gamma) &\le & \sum _{\gamma\in\Gamma_\eps}e^{-\delta(\Gamma) d(\xo,\gamma\xo)}\cdot \frac{\delta(\Gamma)^\alpha h(r_0)}{r_0^\alpha}\cdot e^{\alpha\log d(\xo,\gamma\xo)}\\
&= & \frac{\delta(\Gamma)^\alpha h(r_0)}{r_0^\alpha}\cdot \sum _{\gamma\in\Gamma_\eps}e^{-(\delta(\Gamma)-\alpha) d(\xo,\gamma\xo)}\,.
\ee
 Since $\delta(\Gamma)-\alpha =\frac12\delta(\Gamma)+\frac12 s_\eps >s_\eps$, we conclude that 
$$\sum_{\gamma\in\Gamma_\eps}e^{-b_\gamma} h(b_\gamma)\quad\mbox{converges}.$$ The claim now follows from the fact that $\#\{\gamma\in\Gamma:d(\xo,\gamma\xo)<r_0/\delta(\Gamma)\}$ is finite.\qed

\section{The generalized Patterson-Sullivan construction}\label{genPseries}
% As before $\XX$ is a product of locally compact Hadamard spaces $\XX_1$, $\XX_2$, $\xo=(\xo_1, \xo_2)$ a fixed base point, and 
% $\Gamma\subset\is(\XX_1)\times\is(\XX_2)$ a discrete group which contains two isometries 
% $g=(g_1,g_2)$ and $h=(h_1,h_2)$ \st for $i=1,2$ $g_i$ and $h_i$ are independent rank one elements of $\Gamma_i$. 
According to the statement of Theorem~A, the classical conformal density constructed in the previous section gives measure zero to the set of limit points 
of slope different from $\theta_*$. In order to obtain measures on each $\Gamma$-invariant subset of the limit set we will
%Our aim now is to construct measures for the remaining $\Gamma$-invariant subsets of the limit set which are null sets for the classical conformal densities. 
 use a variation of the classical Patterson-Sullivan construction with more degrees of freedom. The idea is to use a weighted version of the Poincar{\'e} series in order to get the
main contribution from orbit points with slope close to the desired slope $\theta\in (0,\pi/2)$. At this point, properties of the exponent of  growth 
%described in Section~\ref{ExpGrowth} 
and Proposition~\ref{convdiv} will turn out to be  central importance.

Recall that $\bs_\theta$ denotes the directional distance introduced in Section~\ref{prodHadspaces}.
We observe  that for any $b=(b_1, b_2)\in \RR^2$, $\theta\in [0,\pi/2]$  and
$\tau\ge 0$ fixed, the series
$$ P_{\theta}^{s,b,\tau}(x,y)=\sum_{\gamma\in \Gamma} e^{-s
 (b_1 d_1(x_1,\gamma_1 y_1)+b_2 d_2(x_2,\gamma_2 y_2)+\tau(d(x,\gamma y)-\bs_
{\theta}(x,\gamma y )))} $$
possesses a critical exponent which is independent
of $x,y\in\XX$ by the triangle
inequalities  for $d$, $d_1$, $d_2$ and
$\bs_{\theta}$. Notice that for $\tau=0$, this is exactly the series considered by M.~Burger in \cite{MR1230298}; here we will need to take $\tau$ large in order to make the 
contribution of the orbit points with slope far away from $\theta$ small. 
 
For any  $\theta\in [0,\pi/2]$ and $\tau\ge 0$, we 
  define a  {\hl region of convergence}
$$ {\mathcal R}_{\theta}^{\tau}:=\big\{b{=}(b_1, b_2) :
P_{\theta}^{s,b,\tau}(o,o)\ \mbox{has critical exponent}\
s\le1\}\subseteq \RR^2$$
and its boundary
$$ \partial {\mathcal R}_{\theta}^{\tau}:=\big\{b{=}(b_1, b_2) :
P_{\theta}^{s,b,\tau}(o,o)\ \mbox{has critical exponent}\
s\le1\}\subseteq \RR^2.$$
We recall the definition of the distance vector from (\ref{distancevector}). In the sequel we will identify $b=(b_1,b_2)$ with the column vector $b^t$ so that for $x=(x_1,x_2)^t\in\RR^2$ we may write
$$ \langle b,x\rangle =b_1x_1+b_2x_2\,.$$
The region of convergence  possesses the following properties:
\begin{lem}
If $\tau\le \tau'$, then ${\mathcal R}_{\theta}^{\tau}\subseteq {\mathcal R}_{\theta}^{\tau'}$.
\end{lem}
\prf\   Let $\tau\le \tau'$, $b\in {\mathcal
R}_{ \theta}^{\tau}$. Then for any $\gamma\in\Gamma$
$$ e^{-s\big(\langle b, H(\xo,\gamma\xo)\rangle+\tau'(d(o,\gamma o)-\bs_{ \theta}(o,\gamma
    o))\big)}
{\le}\,e^{-s\big(\langle b, H(\xo,\gamma\xo)\rangle+\tau(d(o,\gamma o)-\bs_{ \theta}(o,\gamma
    o))\big)} $$
and therefore 
$P_{ \theta}^{s,b,\tau'}(o,o)\le P_{ \theta}^{s,b,\tau}(o,o)$. Hence
  $P_{ \theta}^{s,b,\tau'}(o,o)$   converges if $s>1$. In particular, $P_{ \theta}^{s,b,\tau'}(o,o)$ possesses a critical
  exponent less than or equal to 1.\qed

\begin{lem}
For any $\tau\ge 0$, the region ${\mathcal R}_{ \theta}^\tau $ is convex.
\end{lem}
\prf\   Let $\tau\ge 0$, $a, b\in {\mathcal R}_{ \theta}^\tau$ and  $t\in[0,1]$.
For $\gamma\in\Gamma$ we abbreviate 
$$\big(ta+(1-t)b\big)_\gamma:=\langle ta+(1-t)b
, H(\xo,\gamma\xo)\rangle +\tau\big(d(o,\gamma
  o)-\bs_{ \theta}(o,\gamma o)\big).$$
Then by H{\"o}lder's inequality 
$$\sum_{\gamma\in \Gamma} e^{-s(t a+(1-t)b)_\gamma}=\sum_{\gamma\in
  \Gamma} e^{-s t a_\gamma} e^{-s (1-t) b_\gamma}\le \Big(\sum_{\gamma\in
    \Gamma}e^{-sa_\gamma}\Big)^t  \Big(\sum_{\gamma\in \Gamma}
  e^{-sb_\gamma}\Big)^{1-t}\,. $$
The latter sum converges if $s>1$, hence $ta+(1-t)b \in {\mathcal
  R}_{ \theta}^\tau$. \qed\\[-1mm]

With the help of Proposition~\ref{convdiv} we can describe the region of convergence more precisely. The following result 
relates the region of convergence ${\mathcal
  R}^\tau_{ \theta}$ to the exponent of growth of slope $\theta$. 
\begin{lem}
Let $\theta\in [0,\pi/2]$   and $\tau\ge 0$. If  $(b_1,b_2)\in {\mathcal R}^\tau_{ \theta}$, then
$\  b_1\cos\theta +b_2\sin\theta  \ge \delta_{ \theta}(\Gamma)\,.$
\end{lem}
\prf\   Recall the definition of $H_\theta$ from~(\ref{thetavector}) and suppose that 
$\langle b,H_\theta\rangle= b_1\cos\theta +b_2\sin\theta < \delta_{ \theta}(\Gamma)$. Then there exists 
$s>1$ \st $s\,\langle b,H_\theta\rangle < \delta_{ \theta}(\Gamma)$. For $H\in\RR^2_{\ge 0}$ we put
$$f(H):=s\big(\langle b, H\rangle+\tau (1-\langle H_\theta,H\rangle)\big)\,.$$
Moreover, the continuous homogeneous function $f:\RR_{\ge 0}^2\to\RR$  satisfies  
$f(H_\theta)=s\langle b, H_\theta\rangle  <
  \delta_{ \theta}(\Gamma)$, so according to  Proposition \ref{convdiv} (a)  applied
to $D=[0,\pi/2]$, the series 
$$\sum_{\gamma\in\Gamma} e^{-f(H(\xo,\gamma\xo))}\qquad\mbox{diverges}\,.$$
Since $f(H(\xo,\gamma\xo))= s\big(b_1 d_1(o_1, \gamma_1 o_1)+b_2 d_2(\xo_2,\gamma_2\xo_2)+\tau(d(\xo,\gamo)-\bs_{ \theta}(o,\gamo)
)\big)$ we obtain a contradiction to $(b_1,b_2)\in {\mathcal R}^\tau_{ \theta}$. \qed\\[-1mm]

Using the above properties of the region of convergence and Patterson's Lemma~\ref{Pt} we are now going to construct \bd ies as defined in
the introduction. Such densities are  a natural generalization of
$\Gamma$-invariant conformal densities if one wants to measure each  $\Gamma$-invariant subset of the geometric limit set.

From here on we fix $\theta\in [0,\pi/2]$ \st $\Lim\cap\rand_\theta\neq \emptyset$, $\tau\ge 0$ and
$b=(b_1,b_2)\in\partial {\mathcal R}^\tau_{ \theta}$. For
$\gamma=(\gamma_1,\gamma_2)\in\Gamma$ we abbreviate  
\begin{equation}\label{bgammadef}
 b_\gamma:=b_1 d_1(o_1,\gamma_1 o_1)+b_2 d_2(\xo_2,\gamma_2\xo_2)+\tau\big(d(o,\gamma o)-\bs_{ \theta}(o,\gamma
  o)\big)\,.
\end{equation}
Let $h$ be a function as in Lemma~\ref{Pt} and recall the definition of the distance vector from (\ref{distancevector}). 
As in Section~\ref{ClasPatSulconst} we will construct a family of
orbital measures on $\overline{\XX}$ in the following way:
If $D$
denotes the unit Dirac point measure, then  for $x \in\XX$ and $s>1$ we put
$$ \mu_x^s:=\frac1{\Summ^s}\sum_{\gamma\in\Gamma}e^{-s\big(\langle b, H(x,\gamma\xo)\rangle+\tau (d(x,\gamma o)-\bs_
{ \theta}(x,\gamma o))\big)}h(b_\gamma)D(\gamma\xo)\,.$$
As in the classical case, these measures are $\Gamma$-equivariant by construction, but they depend on the parameters
$\theta\in [0,\pi/2]$, $\tau\ge 0$ and $b=(b_1,b_2)\in\partial
{\mathcal R}^\tau_{ \theta}$. 

Recall from Section~\ref{ClasPatSulconst} that $(\Cnt^{0}(\ganz),\Vert\cdot\Vert_\infty)$ is the space of real valued continuous functions on $\ganz$ with norm
$\Vert f\Vert_\infty=\max\{|f(x)| :  x\in\ganz\}$,
$f\in \Cnt^{0}(\ganz)$, and $ \rho$ is the pseudo-metric on   
the cone $\MM^+ (\overline\XX)$ of positive finite Borel
measures on $\overline{\XX}$ defined in~(\ref{pseudometricmeasures}).
\begin{lem}\label{equicongen} 
Let $\theta\in [0,\pi/2]$ \st $\Lim\cap\rand_\theta\ne\emptyset$, $\tau\ge 0$ and $b=(b_1,b_2)\in\partial {\mathcal R}^\tau_{ \theta}$. Then the family of maps 
${\mathcal F}(\theta,\tau,b):=\{x\mapsto\mu_x^s :  1<s
\le 2\}$ from $\XX$ to $\MM^+ (\overline\XX)$ is equicontinuous. 
\end{lem}
\prf\  Let $x,y\in\XX$.  For $\gamma=(\gamma_1,\gamma_2)\in\Gamma$ we abbreviate
\begin{multline}\label{kuh}
q_\gamma(y,x):=b_1 
\big(d_1(y_1,\gamma_1\xo_1)-d_1(x_1,\gamma_1\xo_1)\big)+b_2\big(d_2(y_2,\gamma_2\xo_2)-d_2(x_2,\gamma_2\xo_2)\big)\\
+\tau\big(d(y,\gamma
\xo)-d(x,\gamma\xo)
-\bs_{ \theta}(y,\gamma
\xo)+\bs_{ \theta}(x,\gamma\xo)\big)\,,
\end{multline}
put $\Vert
b\Vert_1:= |b_1|+|b_2|$ and estimate
\begin{equation}\label{qest}
\big|q_\gamma(y,x)\big|\le  b_1 d_1(y_1,x_1)+b_2d_2(y_2,x_2)+2\tau
d(y,x)\le  d(x,y)\big(||b||_1+2\tau\big)\,.
\end{equation}
If $s\in(1, 2]$ and $f\in\Cnt^{0}(\ganz)$, the inequality $|1-e^{-t}|\le
e^{|t|}-1$, $t\in\RR$, gives
\begin{align*}
\bigg|\int_{\ganz}
f\;d\mu_x^s - \int_{\ganz} f\;d\mu_y^s\bigg|&\le\frac1{
\Summ^s}\sum_{\gamma\in \Gamma}
e^{-s(b_1 d_1(x_1,\gamma_1\xo_1)+b_2d_2(x_2,\gamma_2\xo_2)+\tau (d(x,\gamma
o)-\bs_{ \theta}(x,\gamma o)))}\\
&\hspace{2.8cm} \cdot h(b_\gamma) |
f(\gamma\xo) |\,\big|1- e^{-s q_\gamma(y,x)}\big|\\
 &\le\frac{\Vert f\Vert_\infty}{\Summ^s}\sum_{\gamma\in 
\Gamma}e^{-s b_\gamma}e^{-sq_\gamma(x,o)} h(b_\gamma)(e^{s
|q_\gamma(y,x)|}-1)\,.\end{align*}
Since $f\in\Cnt^{0}(\ganz)$ was arbitrary, $s\le 2$ and
$\sum_{\gamma\in\Gamma} e^{-sb_\gamma}h(b_\gamma)=\Summ^s$,
we conclude using (\ref{qest}) 
$$
\rho(\mu_x^s,\mu_y^s)\le e^{2d(o,x)(\|b\|_1+2\tau)}\at\big
(e^{2d(x,y)(\|b\|_1+2\tau)}-1\big)\,.$$
This proves that ${\mathcal F}(\theta,\tau,b)$ is equicontinuous.\qed
\begin{lem}
Let $\theta\in [0,\pi/2]$ \st $ \Lim\cap\rand_\theta\neq \emptyset$, $\tau\ge 0$ and $b=(b_1,b_2)\in\partial {\mathcal R}^\tau_{ \theta}$. 
Then for any $x\in\XX$ there exists a sequence 
\hbox{$(s_n)\searrow 1$} \st the measures
 $\mu_{x}^{s_n}\subset \MM^+(\overline\XX)$ converge weakly to a measure $\mu_{x}:=\mu_x( \theta,\tau,b)$ as $n\to\infty$. \end{lem}
\prf\   The compactness of the space $\ganz$ implies that every sequence of measures in $\MM^+(\overline\XX)$ possesses
a weakly convergent subsequence. \qed\\[-1mm]

Hence as in the classical case, the family of maps ${\mathcal F}( \theta,\tau,b)$ is relatively compact in the space of continuous maps $\Cnt(\XX,\MM^+
(\overline\XX))$ endowed with the topology of uniform convergence on compact
sets. 
From the definition of $(\mu_{x}^{s})_{x\in\XX}$ it
follows that every accumulation point
$\mu=\mu( \theta,\tau,b)=(\mu_x)_{x\in\XX}$ of 
${\mathcal F}( \theta,\tau,b)$ 
as $s\searrow 1$  takes its values in  $\MM^+(\rand)$. 

Unfortunately, the families of measures obtained in this way are not very useful, because in general they are not absolutely continuous with respect 
to each other and also depend on the parameter $\tau\ge 0$. However, we still have some freedom in choosing appropriate parameters $b=(b_1,b_2)\in\RR^2$. 

The following  proposition paves the way towards the construction
of orbital measures with support in a single $\Gamma$-invariant subset $ \Lim\cap \rand_\theta\subseteq\rand$.
%%%CONTINUITY POINT????
\begin{prp}\label{key}
Fix $\,\theta\in [0,\pi/2]$ \st $\, \delta_\theta(\Gamma)>0$. %%%IMPLIES $Lim\cap\rand_\theta\ne \emptyset$. 
Then there exist $b=(b_1,b_2)\in\RR^2$ and  $\tau_0\ge 0$ \st for all
$\tau\ge\tau_0$ and for all $ \eps>0$
$$\sum_{\begin{smallmatrix}{\gamma\in\Gamma}\\{|\theta(\xo,\gamma\xo)- \theta|> \eps}\end{smallmatrix}}
e^{-b_\gamma} h(b_\gamma) <\infty\,,$$
where $b_\gamma$ is defined in~(\ref{bgammadef}), and $h$ is a function as in Lemma \ref{Pt}.
\end{prp}
\prf\ Recall the definition of $H_\theta$ from~(\ref{thetavector}). Since the function $\Psi_\Gamma$ defined in (\ref{defPsi}) is concave 
and upper semi-continuous, % is an interior point in the support of $\Summ_\Gamma$
there exists a linear functional $\Phi$ on $\RR^2$ \st
$$\Phi(H_\theta)= \Psi_\Gamma(H_\theta)\quad\an\quad\Phi(H)\ge\Psi_\Gamma(H)
\quad\forall\,H\in\RR_{\ge 0}^2\,.$$ 
We choose $b=(b_1,b_2)\in\RR^2$ \st $\Phi=\langle b, \cdot \rangle$, so
$$ \langle b,H_\theta\rangle =\delta_\theta(\Gamma)\quad\an\qquad \langle b,H_{\hat\theta}\rangle\ge \delta_{\hat\theta}(\Gamma)\quad\forall\, \hat\theta\in [0,\pi/2]\,.$$ 
Since the map $\hat\theta\mapsto \langle b, H_{\hat\theta}\rangle$ is continuous and $\delta_\theta(\Gamma)>0$, there exists $\hat\eps\in( 0,1)$ \st every 
$\hat\theta\in [0,\pi/2]$ with $|\hat\theta-\theta|<\hat\eps$ satisfies $\langle b,H_{\hat\theta}\rangle >0$. 
Notice that it suffices to prove the claim for $\eps<\hat\eps$, because the sum is non-increasing when $\eps$ gets bigger. 

We now put $\Vert b\Vert_1{:=}|b_1|+|b_2|$, fix 
$\tau_0{:=}12\max\{(2\delta(\Gamma){-}\delta_{ \theta}(
\Gamma){+}2
\|b\|_1)/ \hat\eps^2, 2 \}$, and let $ \eps\in (0,\hat\eps)$ be arbitrary.  
By property (ii) of the function $h$ in Lemma \ref{Pt}
there exists $r_0=r_0( \eps)>0$ \st for $r\ge r_0$ and $t>1$ we have  
$h(rt)\le (t)^{ \eps^2} h(r)$. Let $R=R( \eps)> \max\{ \frac{r_0}{2\eps^2}, \frac{r_0}{2\delta(\Gamma)}\}\ $%12 \frac{r_0}{\tau_0 \eps^2},, \frac{r_0}{\delta(\Gamma)}$ 
\st $d(\xo,\gamma\xo)>R$  
implies 
\begin{equation}\label{dbed}
d(\xo,\gamma\xo)\big(\Vert b\Vert_1+2 \eps^2+
  2\delta(\Gamma)\big)<\min\big\{e^{d(\xo,\gamma\xo)},
  e^{\delta(\Gamma)d(\xo,\gamo)}\big\}\,.
\end{equation}
For $\gamma\in\Gamma$ we abbreviate $H_\gamma=H(\xo,\gamma\xo)/d(\xo,\gamma\xo)$, fix $\tau\ge \tau_0$ and  set 
\begin{eqnarray*}
\hat \Gamma & := & \big\{\gamma\in\Gamma : |\theta(\xo,\gamma\xo)-\theta|>\frac{\hat\eps}2\},\\
\Gamma'&:= &\{\gamma\in\Gamma : \frac{\eps}2 < |\theta(\xo,\gamma\xo)-\theta|<\hat\eps\};
\end{eqnarray*}
our goal is to show that $\sum_{\gamma\in \hat \Gamma} e^{-b_\gamma} h(b_\gamma)<\infty$ and $\sum_{\gamma\in\Gamma'}   e^{-b_\gamma} h(b_\gamma)<\infty$.

First let $\gamma\in \hat\Gamma$. Using the Cauchy--Schwarz
inequality, $\Vert H_\gamma - H_\theta\Vert\le 2$, and the condition $\,\langle b,H_\theta \rangle = 
\delta_{ \theta}(\Gamma)$,  we have
\begin{eqnarray*}
\langle b,H_\gamma\rangle&=& \langle b, H_\theta\rangle +  
\langle b, H_\gamma- H_\theta\rangle \ge \delta_{ \theta}(\Gamma)-\sqrt{b_1^2+b_2^2}\cdot  \Vert H_\gamma- H_\theta  \Vert \ge \delta_{ \theta}(\Gamma)-2\|b\|_1\,.
\end{eqnarray*}
Moreover, the estimate $\ \cos t<1-t^2/3\ $ for $\, t\in\RR$,  and the fact that \\$ d(\xo,\gamma\xo)-\bs_\theta(\xo,\gamma\xo)>d(\xo,\gamma\xo)\big(1-\cos\frac{\hat\eps}{2}\big)$ and    
$\tau\ge 12(2\delta(\Gamma)-\delta_{ \theta}(\Gamma)+2
\|b\|_1)/ \hat\eps^2$ imply 
\begin{eqnarray*}
b_\gamma&= & \langle b,H(\xo,\gamma\xo) \rangle +\tau\big(d(\xo,\gamma\xo)-\bs_\theta(\xo,\gamma\xo)\big)\\
& > &  d(\xo,\gamo)\left(\delta_{ \theta}(\Gamma)-2\|b\|_1+
\frac{\tau \hat\eps^2}{12}\right)\\
\ge  2\delta(\Gamma) d(\xo,\gamo)\,.
\end{eqnarray*}
Hence if $\gamma\in\hat\Gamma$ satisfies $d(\xo,\gamma\xo)>R$, then by choice of $R>\frac{r_0}{2\delta(\Gamma)}$ we have $b_\gamma > r_0$, hence 
\begin{equation}\label{habsch}
 h(b_\gamma)=h\left(\frac{b_\gamma}{r_0} r_0\right)\le
\left(\frac{b_\gamma}{r_0}\right)^{ \eps^2}h(r_0) =
\frac{h(r_0)}{(r_0)^{ \eps^2}}e^{ \eps^2\log b_\gamma}
\,.
\end{equation}
Since the function $g(t):= t- \eps^2\log t$, $t>0$,  
is monotone increasing, we conclude that for $\gamma\in \hat\Gamma$ with $d(\xo,\gamma\xo)>R$  
\begin{eqnarray*}
 g(b_\gamma)  & >  & g\big(2\delta(\Gamma)d(\xo,\gamma\xo)\big)
= 2\delta(\Gamma)d(\xo,\gamma\xo) -\eps^2\log \big(\underbrace{2\delta(\Gamma)d(\xo,\gamma\xo)}_{<e^{\delta(\Gamma) d(\xo,\gamma\xo)}}\big)\\
&>& d(\xo,\gamma\xo)\big(2\delta(\Gamma)-\eps^2\delta(\Gamma)\big)\,,
\end{eqnarray*}
where we used (\ref{dbed}) in the last step.  Summarizing, we estimate 
\begin{eqnarray*}
 \sum_{\begin{smallmatrix}{\gamma\in\hat\Gamma}\\
{d(\xo,\gamo)>R}\end{smallmatrix}}\!\!e^{-b_\gamma}h(b_
\gamma) &\le & \frac{h(r_0)}{(r_0)^{ \eps^2}}\sum_{
\begin{smallmatrix}{\gamma\in\hat\Gamma}\\{d(\xo,\gamo)>R}
\end{smallmatrix}}\!\!e^{-b_\gamma+ \eps^2\log b_\gamma}
=\frac{h(r_0)}{(r_0)^{ \eps^
2}}\sum_{\begin{smallmatrix}{\gamma\in\hat\Gamma}\\{d(\xo,
\gamo)>R}\end{smallmatrix}}
\!\!e^{-g(b_\gamma)}\\
& < & 
\frac{h(r_0)}{(r_0)^{ \eps^2}}\sum_{\gamma\in\Gamma}e^{-
(2- \eps^2)\delta(\Gamma)  d(\xo,\gamo)}\,,
\end{eqnarray*}
which converges since $ \eps^2\le  \hat\eps^2 <1$.

Next let $\gamma\in\Gamma'$. As above, $|\theta(\xo,\gamma\xo)-\theta|>\eps/2$ and $\tau\ge 24$ imply
$$b_\gamma >\langle b,H(\xo,\gamma\xo)\rangle +\tau d(\xo,\gamma\xo)\frac{\eps^2}{12}>d(\xo,\gamma\xo)\big(\langle b,H_\gamma\rangle +2\eps^2\big)\,,$$
hence % if $d(\xo,\gamma\xo)>R$, then by (\ref{dbed})
$$  g(b_\gamma)> d(\xo,\gamma\xo)\big( \langle b,H_\gamma\rangle +2\eps^2\big)-\eps^2\log\big( d(\xo,\gamma\xo)( \langle b,H_\gamma\rangle +2\eps^2)\big) \,.$$
By choice of $\hat\eps$ we have $\langle b, H_\gamma\rangle>0$, so $d(\xo,\gamma\xo)>R>\frac{r_0}{2\eps^2}$ yields $ b_\gamma>r_0$. Moreover, (\ref{dbed}) implies 
$d(\xo,\gamma\xo)\big(\langle b,H_\gamma\rangle +2\eps^2\big)<e^{d(\xo,\gamma\xo)}$, which gives
$$g(b_\gamma)>d(\xo,\gamma\xo)\big(\langle b,H_\gamma\rangle +\eps^2\big)\,.$$
Next we consider the continuous homogeneous function \[f: \RR^2_{\ge 0}\to\RR,\ H\mapsto \langle b,
H \rangle + \eps^2\Vert H\Vert.\]  
Using inequality (\ref{habsch}) we estimate
\begin{align*} 
\sum_{\begin{smallmatrix}\gamma\in\Gamma'\\
d(\xo,\gamo)>R\end{smallmatrix}}\!\!e^{-b_\gamma}h(b_
\gamma)&
\le\frac{h(r_0)}{(r_0)^{ \eps^2}}\sum_{\begin{smallmatrix}{\gamma\in\Gamma'}
\\
{d(\xo,\gamo)>R}\end{smallmatrix}}
\!\!e^{-g(b_\gamma)}<\frac{h(r_0)}{(r_0)^{ \eps^2}}\sum_{
\begin{smallmatrix}{\gamma
\in\Gamma'}\\{d(\xo,\gamo)>R}\end{smallmatrix}}
\!\!e^{-f(H(\xo,\gamma\xo)) }\\
&<\frac{h(r_0)}{(r_0)^{ \eps
^2}}\sum_{\gamma\in\Gamma'}\!\!e^{-f(H(\xo,\gamma\xo))}\,. 
\end{align*}
This sum converges by Proposition \ref{convdiv} (b)
applied to \[D':=\{\hat\theta\in [0,\pi/2]: \eps/2 <|\hat\theta -\theta|<\hat\eps\},\] because
the continuous homogeneous function
$f$ satisfies $f(H_{\hat\theta})>\delta_{\hat\theta}(\Gamma)$ for all
$\hat\theta\in \overline{D'}$.

The claim now follows from the fact that the number of elements $\gamma\in\Gamma$ with $d(\xo,\gamma\xo)\le R$ is finite. \qed\\[-1mm]

Notice that the above proposition provides $b=(b_1,b_2)\in\RR^2$ even if $\theta=0$ or $\theta=\pi/2$. However, in general we do not have $b_2=0$ 
if $\theta=0$, or $b_1=0$ if $\theta=\pi/2$. 
Hence for the construction of \bd ies we have to restrict ourselves to $\theta\in (0,\pi/2)$. \\

We therefore fix $\theta\in (0,\pi/2)$ \st $\delta_\theta(\Gamma)>0$, hence in particular we have $\Lim\cap\rand_\theta\neq\emptyset$.  
Using the previous proposition we are finally able to construct a \bd y according to Definition~\ref{bdensi}.

To that end we choose  
$b=(b_1,b_2)\in\RR^2$ and $\tau_0>0$ according to Proposition~\ref{key}.
For fixed  $\tau>\tau_0$ we consider  the corres\-pon\-ding family ${\mathcal F}(\theta,\tau,b)$
as in {\rm Lemma} \ref{equicongen}.
Then by the Theorem of Arzel{\`a}-Ascoli,
\hbox{${\mathcal F}(\theta,\tau,b)$} is relatively
compact in the space of continuous maps $\Cnt(\XX,\MM^+
(\overline\XX))$ endowed with the topology of uniform convergence on compact
sets. The following
proposition characterizes the  possible accumulation points.
\begin{prp} Every accumulation point $\mu=\mu(\theta,\tau,b)$ of the
family ${\mathcal F}={\mathcal F}(\theta,\tau,b)$ in  $\Cnt(\XX,\MM^+
(\overline\XX))$  is a \bd y.
\end{prp}
\prf\ Let $(\mu_x)_{x\in\XX}$ be an accumulation point of ${\mathcal
 F}$. By construction, the measures $\mu_x$, $x\in\XX$, are
 $\Gamma$-equivariant and supported on the limit set $\Lim$. 
Proposition~\ref{key} further implies
$\supp(\mu_\xo)\subseteq
 \Lim\cap \rand_\theta$. 
It therefore suffices to prove 
$$\frac{d\mu_x}{d\mu_\xo}(\tilde\eta)=e^{b_1\bs_{\eta_1}(\xo_1,x_1) + b_2\bs_{
\eta_2}(\xo_2,x_2)}\quad\mbox{for any}\ x\in\XX\,,\ \tilde\eta=(\eta_1,\eta_2,\theta)
\in\supp(\mu_\xo)\,.$$
Notice that if $(y_n)=\big((y_{n,1},y_{n,2})\big)\subset \XX$ is a sequence converging to 
a point $\tilde\eta=(\eta_1,\eta_2,\theta) \in \rand_\theta$,  then for $i\in\{1,2\}$ 
$d_i(x_i,y_{n,i})-d_i(\cdot, y_{n,i})\to \bs_{\eta_i}(x_i,\cdot)$ uniformly on compact sets in $\XX_i$, and 
$d(x,y_n)-d(\cdot , y_n)\to \bs_{\tilde\eta}(x,\cdot)$ uniformly on compact sets in $\XX$.  Using the definition and properties of 
Busemann functions and the directional distance $\bs_\theta$ we conclude that  for any constant $c\ge 0$ and $\eps>0$ arbitrary,
there exist $R>0$ and $\rho>0$ with the following properties: If $\gamma=(\gamma_1,
\gamma_2) \in\Gamma$ satisfies  $d(o,\gamma o)>R$ and  $|\theta(\xo,\gamma\xo)-\theta|<\rho$, then for
all $x=(x_1,x_2) \in\XX$ with $d(\xo,x)\le c$
\begin{equation}\label{unco} 
\big|d(o,\gamma\xo){-}d(x,\gamma\xo){-}\bs_{\theta}(o,\gamma\xo)+\bs_\theta(x,\gamma\xo)\big| < \eps\,.
\end{equation}
Let $\eps>0$ arbitrary, fix $x\in\XX$, put $c:=d(o,x)$ and choose $R>0$ 
and $\rho >0$ as above. 
We set \vspace{-4mm}
\begin{align*}
%&\Gamma_1:=\big\{\gamma\in\Gamma : d(o,\gamma\xo)\le R
%\big\}\,,\\
&\hat\Gamma:=\big\{\gamma\in\Gamma : \big| \theta( o, \gamma \xo)-\theta\big|>\rho/2\big\}\,,\\
&\Gamma':=\big\{\gamma\in\Gamma : d(o,\gamma\xo)>
R\,\an\  \big| \theta(o,\gamma\xo)-\theta\big|<\rho\big\}\,.
\end{align*}
Here, for $\gamma\in\Gamma$ we abbreviate 
$$q_\gamma(o,x):= d(\xo,\gamma\xo)-d(x,\gamma\xo)-\bs_\theta(o,\gamma\xo)+\bs_\theta(x,\gamma\xo)$$
and recall that for $\gamma\in\Gamma'$ we have $|q_\gamma(\xo,x)|<\eps$. 
Then for any $f\in \Cnt^0(\ganz)$, $s\in (1,2]$, we have  
\begin{eqnarray*}
\bigg|\int_{\ganz} f(z)d\mu_\xo^s(z)& - & \int_{\ganz}f(
z)e^{-b_1(d_1(\xo_1,z_1)-d_1(x_1,z_1))-b_2 (d_2(o_2,z_2)-d_2(x_2,z_2))}d\mu_x^s(z)
\bigg|\\
& \le & \frac1{\Summ^s}\sum_{\gamma\in\Gamma}|f(\gamma\xo)|
\cdot e^{-sb_\gamma}h(b_\gamma)\cdot |1-e^{(s-1)\langle b, H(\xo,\gamma\xo)-H(x,\gamma\xo)\rangle+s\tau q_\gamma(o,x)}\big|\,.
\end{eqnarray*}
The triangle inequality and the estimate $|q_\gamma(\xo,x)|\le 2 d(\xo,x)$ imply
that for any $\gamma\in\Gamma$
$$ \big|(s-1)\langle b, H(\xo,\gamma\xo)-H(x,\gamma\xo)\rangle+s\tau q_\gamma(o,x)\big| \le (s-1)\Vert b\Vert_1 d(o,x)+2s\tau d(o,x)\,.$$
This proves
that for $x\in\XX$ with $d(x,\xo)\le c$, and $s\le 2$, the term 
\begin{eqnarray*}
\big|1-e^{(s-1)\langle b, H(\xo,\gamma\xo)-H(x,\gamma\xo)\rangle+sq_\gamma(o,x)}\big|
&\le & e^{|(s-1)\langle b, H(\xo,\gamma\xo)-H(x,\gamma\xo)\rangle+sq_
\gamma(o,x)|}-1\\
&\le& e^{c(\|b\|_1+4\tau)}-1
\end{eqnarray*}
is bounded above by a constant $A=A(c,b,\tau)$. 
If $\gamma\in\Gamma'$ then $|q_\gamma(\xo,x)|<\eps$, hence
\begin{equation}\label{limes} \lim_{s\searrow 1}
\big|1-e^{(s-1)\langle b, H(\xo,x)\rangle+sq_
\gamma(o,x)}\big|\le \lim_{s\searrow 1} e^{(s-1)\Vert b\Vert_1 c +2\tau\eps}
  -1=  e^{2\tau\eps}-1.\end{equation}
We conclude 
\goodbreak\noindent
\begin{multline*}
\bigg|\int_{\ganz}\!f(z)d\mu_\xo^s(z){-}\int_{\ganz}
\!f(z)e^{-b_1(d_1(\xo_1,z_1)-d_1(x_1,z_1))-b_2(d_{2}(o_2,z_2)-d_2(x_2,z_2))}d\mu_x^s(z)
\bigg|\\
{\le}  \frac{\Vert f\Vert_\infty}{\Summ^s}
\Big(\!\sum_{\begin{smallmatrix}{\gamma\in\Gamma}\\{d(\xo,\gamma\xo)\le R}\end{smallmatrix}}
e^{-sb_\gamma}h(b_\gamma) 
+
\sum_{\gamma\in\hat\Gamma} e^{-sb_\gamma}h(b_\gamma) A\\*
+\sum_{\gamma\in\Gamma'}
e^{-sb_\gamma}h(b_\gamma)
\big|1-e^{(s-1)\langle b,H(\xo,\gamma\xo)-H(x,\gamma\xo)\rangle+s\tau q_
\gamma(o,x)}\big|\Big)\,.
\end{multline*}
Now the first term tends to zero as $s\searrow 1$ since the summation is over a finite number of $\gamma\in\Gamma$. By Proposition~\ref{key}, 
$\sum_{\gamma\in\hat\Gamma} e^{-b_\gamma}h(b_\gamma)$
converges, hence the second term tends to zero as
$s\searrow 1$.
In the last term, we have $\sum_{\gamma\in\Gamma'}
e^{-sb_\gamma}h(b_\gamma)\le \Summ^s$ for any $s>1$, and therefore by (\ref{limes})
\begin{multline*}
\lim_{s\searrow 1}\bigg|\int_{\ganz}f(z)d\mu_\xo^s(z
)-\int_{\ganz}f(z)e^{-b_1(d_1(\xo_1,z_1)-d_1(x_1,z_1))-b_2(d_2(o_2,z_2)-d_2(x_2,z_2))}
d\mu_x^s(z)\bigg|\\
\le\Vert f\Vert_\infty\big(  e^{2\tau\eps}-1\big)\,.
\end{multline*}
The claim follows
taking the limit as $\eps\searrow 0$. \qed

So in particular we can construct for any $\theta\in (0,\pi/2)$ with $\delta_\theta(\Gamma)>0$ a $(b,\theta)$-density with appropriate parameters $b=(b_1,b_2)$. This proves Theorem~B
from the introduction.

\section{Properties of  
{\boldmath{$(b,\theta)$}}-densities}\label{Propbdies}

In this section we will study properties of \bd ies using the shadow lemma Theorem~\ref{shadowlemma}. If not otherwise specified we allow $\theta\in [0,\pi/2]$.
% $\Gamma\subset\is(\XX_1)\times\is(\XX_2)$  is again a discrete group acting on a product of locally compact Hadamard spaces $\XX_1$, $\XX_2$, which contains two isometries 
% $g=(g_1,g_2)$ and $h=(h_1,h_2)$ \st $g_i$ and $h_i$ are independent rank one elements of $\Gamma_i$  for $i=1,2$. For the sequel we allow $\theta\in [0,\pi/2]$ 
% if not otherwise specified. We will study properties of \bd ies using the shadow lemma Theorem~\ref{shadowlemma}. 

\begin{lem}
Let $\mu$ be a \bd y, and $x\in\XX$.  If $\,\tilde U\subset\rand$ is an open neighborhood of a limit point $\tilde\xi\in\rand_\theta$, then  
$\mu_x(\tilde U) >0$. 
\end{lem}
\prf\ Let $\tilde U \subset\rand$ be an  open neighborhood of a limit point $\tilde\xi\in\rand_\theta$ \st $\mu_x(\tilde U)=0$. If $U:=\tilde U\cap\rand_\theta$, then by compactness 
and minimality of $\Lim\cap\rand_\theta$ (see Theorem~A in \cite{MR2629900}) there exists a finite set $\Lambda\subset\Gamma$ \st 
$$\Lim\cap\rand_\theta  \subseteq \bigcup_{\gamma\in\Lambda} \gamma U\,.$$
Moreover, by $\Gamma$-equivariance
$$ \mu_x(\Lim\cap\rand_\theta)\le\sum_{\gamma\in\Lambda}\mu_x(\gamma U)=\sum_{\gamma\in\Lambda}\mu_{\gamma^{-1}x}(U)\le \sum_{\gamma\in\Lambda}\mu_{\gamma^{-1}x}(\tilde U)=0\,,$$
since $\mu_{\gamma^{-1}x}$, $\gamma\in\Lambda$, is absolutely continuous with respect to $\mu_x$.\qed\\[-1mm]

Recall the definition of the distance vector (\ref{distancevector}) from Section~\ref{prodHadspaces}. 
\begin{thr}[\,Shadow lemma]\label{shadowlemma} 
Let $\mu$ be a \bd y. Then there exists a constant $c_0>0$ \st
for any $c>c_0$ %and for every $\Gamma$-invariant measurable set {\gamma\in\Gamma}\\$F\subseteq\rand$ with $\mu_\xo(F)>0$ 
there exists a constant  $D(c)>1$ with the property 
$$ \frac1{D(c)}e^{-\langle b, H(\xo,\gamma\xo)\rangle} \le \mu_{o
}\big(\Sh(o:B_{\gamma o}(c))\big)\le D(c)e^{-\langle b, H(\xo,\gamma\xo)\rangle}.$$
\end{thr}   
\prf\   Let $U_1\subset\rand_1$, $U_2\subset\rand_2$ be neighborhoods of $h_1^+$, $h_2^+$, $\Lambda\subset\Gamma$ finite, and $c_0>0$ as in 
Proposition~\ref{largeshadows}. For $\theta\in [0,\pi/2]$ and 
$\alpha=(\alpha_1,\alpha_2)\in\Lambda$ the  sets 
$$U_\alpha:=  \Bigg\{\begin{array}{lcl}\ \{(\eta_1,\eta_2,\theta)\in\rand_\theta :\eta_1\in \alpha_1 U_1, \eta_2\in\alpha_2 U_2\}  & \mbox{if} & \theta\in (0,\pi/2)\,\\
 \quad  \alpha_1 U_1 \subset \rand_1\cong (\rand)_0  & \mbox{if} & \theta=0\,\\
\quad \alpha_2 U_2\subset\rand_2\cong (\rand)_{\pi/2} & \mbox{if} & \theta=\pi/2\end{array}$$ 
are relatively open neighborhoods of a limit point in $\rand_\theta$, so by the previous lemma   
$$ q:=\min\{ \mu_\xo(U_\alpha): \alpha\in\Lambda\}\,$$
is strictly positive. Moreover, if $c\ge c_0$ and $\gamma\in\Gamma$ \st $d(o,\gamo)>c$ then by Proposition~\ref{largeshadows} there exists $\alpha\in\Lambda$ \st $U_\alpha\subseteq \Sh(\gamma^{-1}\xo:B_\xo(c))$. Hence for $c\ge c_0$ and
$\gamma\in\Gamma$ with $d(o,\gamo)>c$ we have
\begin{equation}\label{mush}
\mu_{o}(\rand)\ge\mu_{o}\big(\Sh(\gamma^{-1}o :B_{o}(c))
\big)\ge q>0\,.\end{equation}
Put $S_\gamma:=\Sh(\xo:B_{\gamma o}(c))$ and recall the definition of the $b$-Busemann function (\ref{buscomb}). The properties (ii) and (iii) of a \bd y %and the $\Gamma$-invariance of $F$ 
imply
\begin{align*}
\mu_{o}\big(\Sh(\gamma^{-1}\xo : B_{o}(c))\big)&=\mu_{o}(\gamma^{-1}
S_\gamma) =\mu_\gamo(S_\gamma)= \int_{S_\gamma}d\mu_{\gamo}(\tilde\eta)\\
&=\int_{\smSh(o:B_{\gamo}(c))}e^{\bs^b_{\tilde\eta}(\xo,\gamma\xo)} d\mu_{o}(\tilde\eta)\,.
\end{align*}
%Notice that thanks to our additional conditions for \bd ies in the case $\theta=0$ or $\theta=\pi/2$ the expression in the exponent of the last term is well-defined, even if $\bs_{\eta_1}$ 
%or $\bs_{\eta_2}$ is not.  
By {\rm Lemma}~\ref{esti},
$$ 
e^{-2c}e^{\langle b, H(\xo,\gamma\xo)\rangle }\mu_{o}(S_\gamma)
< \mu_{o}\big(\Sh(\gamma^{-1}o:B_{o}(c))\big)
\le  e^{\langle b, H(\xo,\gamma\xo)\rangle} \mu_{o}(S_\gamma)\,,
$$
so equation (\ref{mush}) allows us to conclude
$$
e^{-\langle b, H(\xo,\gamma\xo)\rangle} q
%\left(\mu_{o}(F)-q\right)
\le \mu_{o}(S_\gamma)\le e^{-\langle b, H(\xo,\gamma\xo)\rangle} e^{2c}\at\mu_{o}(\rand)\,.\eqno{\scriptstyle\square}
$$

The following applications of Theorem~\ref{shadowlemma} yield relations between the  exponent of growth of a given slope 
$\theta\in [0,\pi/2]$ and the parameters of a \bd y. Recall the definition of $H_\theta$ from (\ref{thetavector}).

\begin{thr}\label{deltaklein}
If for $\theta\in (0,\pi/2)$ a $\Gamma$-invariant
$(b,\theta)$-density exists, then
$$\delta_{\theta} (\Gamma)\le \langle b, H_\theta\rangle.$$
\end{thr}
\prf\   Suppose $\mu$ is a \bd y. Let
$c>c_0+1$, where $c_0>0$ is as in {\rm Theorem}~\ref{shadowlemma}, $\eps>0$ and $R>3c_0$  arbitrary. Let $\tilde\eta=(\eta_1,\eta_2,\theta)\in\supp(\mu_\xo)$.
We only need $N(\eps) R$ balls of radius 1 in
$\XX$ to cover the set 
$$\{\big(\sigma_{\xo_1,\eta_1}(t\cos\hat\theta),\sigma_{\xo_2,\eta_2}(t\sin\hat\theta)\big)\in\XX: R-1\le t <R,\, |\hat\theta -\theta|<\eps\}\,,$$
and
$N(\eps)$ is independent of $R$. Since $\Gamma$ is discrete, a
$2c$-neighborhood of any of these balls contains a uniformly bounded number
$M_c$ of ele\-ments of $\Gamma\cdot\xo$.

The compactness of $\rand_\theta$ implies the
existence of  a constant $A>0$ \st every point in $\rand_\theta$ is contained in at
most $A M_c N(\eps) R$ sets 
$\Sh(\xo: B_\gamo(c))$, $\gamma\in \Gamma'$, where
$\Gamma':=\{\gamma\in\Gamma :|\theta(\xo,\gamma\xo)-\theta|<\eps,\,R-1\le
    d(\xo,\gamma\xo)<R\}$.
Therefore 
$$\sum_{\gamma\in\Gamma'}\mu_\xo\big(\Sh(\xo: B_
\gamo(c))\big) \le A M_c 
    N(\eps) R \mu_\xo(\rand_\theta)=A M_c 
    N(\eps) R \mu_\xo(\rand)\,.$$
Furthermore, if $\, \gamma\in\Gamma'$ then $H_\gamma:=H(\xo,\gamma\xo)/d(\xo,\gamma\xo)\, $ satisfies $ \Vert H_\gamma -H_\theta\Vert\le\eps$. 
Writing $\Vert b\Vert_1:=|b_1|+|b_2|$, 
using the Cauchy--Schwarz inequality and $\sqrt{b_1^2+b_2^2}\le \Vert b\Vert_1$ we obtain for $\gamma=(\gamma_1,\gamma_2)\in\Gamma'$
$$
\langle b,H_\gamma\rangle =\langle b,H_\theta\rangle +\langle b, H_\gamma-H_\theta\rangle\le \langle b,H_\theta\rangle   +\Vert b\Vert_1 \eps\,.$$
Using {\rm Theorem}~\ref{shadowlemma} and
$$\Delta N_\theta^\eps(\xo,\xo;R):=\#\{ \gamma\in\Gamma\;:\, R-1\le d(\xo,\gamma \xo)<R\,,\ |\theta(\xo,\gamma\xo)-\theta|<\eps\}\,,\quad R\gg 1 \,,$$ 
we conclude
\begin{multline*}
\Delta N_{\theta}^\eps(o,o;R)\frac1{D(c)}e^{-\langle b,H_\theta\rangle R}{\le}\sum_{\gamma
\in\Gamma'}\frac1{D(c)}
    e^{-\langle b, H(\xo,\gamma\xo)\rangle +\eps \Vert b\Vert_1  d(\xo,\gamo)}\\
 \le  e^{\eps \Vert b\Vert_1
    R}\sum_{\gamma\in\Gamma'}\mu_\xo\big(\Sh(\xo:B
_\gamo(c))\big) 
\le e^{\eps \Vert b\Vert_1
    R} A M_c 
    N(\eps) R \mu_\xo(\rand)\,.
\end{multline*}
Hence
\begin{eqnarray*}
\delta_{\theta}^\eps(o,o)& \le & \limsup_{R\to\infty}
\frac1{R}\log \bigg( D(c) A M_c 
N(\eps) \mu_\xo(\rand) R\cdot\exp\big(\langle b,H_\theta\rangle R+\eps \Vert b\Vert_1
    R\big)\bigg)\\
&=& \langle b,H_\theta\rangle + \eps \Vert b\Vert_1
\end{eqnarray*}
and the claim follows as $\eps\searrow 0$.\qed\\[-1mm]

Unfortunately, the proof of the above proposition does not work for $\theta\in \{0,\pi/2\}$. So we do not know whether the estimate holds for $\delta_{0}(\Gamma)$ and $\delta_{\pi/2}(\Gamma)$.

We next recall the  notion of the  radial limit set from Definition \ref{raddef} of the introduction. If $\theta\in [0,\pi/2]$, then the radial limit set in $\rand_\theta$ is given by
\begin{equation}\label{radlimthetacover}
 \radlim\cap\rand_\theta=\bigcup_{c>0}\bigcap_{R>c}\, \bigcap_{\eps>0}
\bigcup_{\begin{smallmatrix}{\scriptscriptstyle\gamma\in
\Gamma}\\
{\scriptscriptstyle d(o,\gamma
\xo)>R}\\{\scriptscriptstyle |\theta(\xo,\gamma\xo)-\theta|<
\eps}\end{smallmatrix}}
\Sh(\xo: B_{\gamma\xo}(c))\cap\rand_\theta\,.  
\end{equation}

Together with the previous theorem the following implies that if a \bd y gives positive measure to
the regular radial limit set, then the exponent of growth of $\Gamma$ of slope $\theta$ is completely determined by its parameters.
\begin{thr}\label{posmeasrad}
If $\theta\in [0,\pi/2]$ and a \bd y gives positive measure to $\radlim$, then 
$ \delta_{\theta}(\Gamma) \ge  \langle b, H_\theta\rangle$.
\end{thr}
\prf\   Suppose  $\mu$ is a \bd y  \st $\mu_\xo(\radlim)>0$.
Let $c>c_0$  with $c_0>0$ as in Theorem~\ref{shadowlemma}.
Let $\eps>0$ and $R>c$ arbitrary,  and set 
\[\Gamma':=\{\gamma\in\Gamma :  d(o,\gamma\xo)>R,\ | \theta(\xo,\gamma\xo)-\theta|<\eps\}.\] 
Then by (\ref{radlimthetacover})
$$\radlim\cap\rand_\theta \subseteq \bigcup_{\gamma\in\Gamma'}
\Sh(\xo: B_\gamo (c))\cap\rand_\theta\,,$$
and  we estimate 
$$ 0<\mu_{\xo}(\radlim)=\mu_\xo(\radlim\cap \rand_\theta)\le\sum_{\gamma\in\Gamma'}\mu_\xo\big(\Sh(\xo: B_\gamo(c))\big) \le D(c)\sum_{\gamma\in\Gamma'}e^{-\langle b, H(\xo,\gamma\xo)\rangle}.$$
This  implies that for any
$\eps>0$ the tail of the series
$$
\sum_{\begin{smallmatrix}{\gamma\in
\Gamma}\\{|\theta(\xo,\gamo)-\theta|<\eps}\end
{smallmatrix}}e^{-\langle b, H(\xo,\gamma\xo)\rangle}$$ 
does not tend to zero. Therefore the sum above diverges,
and by {\rm Proposition} \ref{convdiv} (b) there exists
$\hat\theta\in [0,\pi/2]$, $|\hat\theta-\theta|\le\eps$ \st
$$\langle b,H_{\hat\theta}\rangle \le
\delta_{\hat\theta}(\Gamma)\,.$$ 
Taking the limit as $\eps\searrow 0$, we conclude 
$ \langle b,H_\theta\rangle \le\delta_{\theta}
(\Gamma)$. \qed\\[-1mm]

Recall the definition of the $b$-Busemann function (\ref{buscomb}) from Section~\ref{prodHadspaces}. The following two lemmata hold for any 
 $\theta\in [0,\pi/2]$ 
and will be important for the proof of Theorem~\ref{radpoint}.% B.
\begin{lem}\label{busnull}
Let $\mu$ be a $( b,\theta)$-density. If $\tilde\eta\in \rand_\theta$ is a point mass for $\mu$,  and 
$\Gamma_{\tilde\eta}$ its stabilizer, then for 
any $\gamma\in\Gamma_{\tilde\eta}$ and $x\in \XX\;$  we
have $$ \bs_{\tilde\eta}^b(x,\gamma x)= b_1\bs_{\eta_1}(x_1,\gamma_1x_1)+b_2 \bs_{\eta_2}(x_2,\gamma_2 x_2)=0\,.$$
In particular, if $\gamma,\hat\gamma\in\Gamma$ are representatives of the same
coset in $\Gamma/\Gamma_{\tilde\eta}$, then % and $\eta:=\pi^F(\tilde\eta)$, then 
 $$ \bs^b_{\tilde\eta}(x,\gamma^{-1}x)= \bs^b_{\tilde\eta}(x,\hat\gamma^{-1}x) \,.$$
\end{lem}
\prf\  If $\gamma\in \Gamma_{\tilde\eta}$, then for $x\in\XX$ we have by $\Gamma$-equivariance 
$$ \mu_x(\tilde \eta)=\mu_x(\gamma^{-1}\tilde\eta)=\mu_{\gamma x}(\tilde\eta)\,.$$
From the assumption that $\tilde\eta$ is a point mass and 
property (iii) in Definition~\ref{bdensi} we get
$$ 1=\frac{\mu_{\gamma x}(\tilde\eta)}{\mu_{x}(\tilde\eta)}=e^{\bs^b_{\tilde\eta}(x,\gamma x)}\,,$$
hence  $\bs^b_{\tilde\eta}(x,\gamma x) =0$ for any $\ging_{\tilde\eta}$.
 
Let $\gamma,\hat \gamma\in\Gamma$ \st
$\gamma\Gamma_{\tilde\eta}=\hat\gamma\Gamma_{\tilde\eta}\in \Gamma/\Gamma_{\tilde\eta}$. Then
$\hat\gamma^{-1}\gamma\in\Gamma_{\tilde\eta}$ and we obtain from 
the above, using the cocycle identity for the Busemann functions $\bs_{\eta_1}$, $\bs_{\eta_2}$,
\begin{align}
\bs^b_{\tilde\eta}(x,\gamma^{-1}x)&= \bs^b_{\tilde\eta}(x,\gamma^{-1}  x) +\bs^b_{\tilde\eta}(\gamma^{-1}x,\hat\gamma^{-1}\gamma \gamma
^{-1}x)\nonumber\\
&=\bs^b_{\tilde\eta}(x,\hat\gamma^{-1}\gamma \gamma^{-1}  x) = \bs^b_{\tilde\eta}(x,\hat\gamma^{-1}
x)\,.\tag*{$\scriptstyle\square$}
\end{align}

\begin{lem}\label{sumconv}
If $\tilde\eta\in
\rand_\theta$ is a point mass for a $( b,\theta)$-density $\mu$, then the sum
$$ \sum e^{\bs^b_{\tilde\eta}(\xo,\gamma^{-1} o)}$$
taken over a system of coset representatives of $\Gamma/\Gamma_{\tilde\eta}$
 converges. \end{lem}
\prf\   If $\gamma$ and $\hat\gamma$ are representatives of different cosets in
$\Gamma/\Gamma_{\tilde\eta}$, then
$\gamma\tilde\eta\neq \hat\gamma\tilde\eta$  and so, by $\Gamma$-equivariance, 
the sum
$\sum\mu_{\gamma^{-1}\xo}(\tilde\eta)=\sum\mu_\xo(\gamma\tilde\eta)$
over a system of coset representatives of $\Gamma/\Gamma_{\tilde\eta}$ is bounded above by
$\mu_o(\rand)$. 
By property (iii) in Definition \ref{bdensi} and the
assumption that
$\tilde\eta$ is a point mass we conclude that the sum
$$ \sum  e^{\bs^b_{\tilde\eta}(\xo,\gamma^{-1} o)} =
  \sum \frac{\mu_{\gamma^{-1}\xo} (\tilde\eta)}{\mu_{\xo}(\tilde\eta)}  =
\frac1{\mu_\xo(\tilde\eta)} \sum  \mu_{\gamma^{-1}\xo} (\tilde\eta)     $$
over a system of coset representatives of $\Gamma/\Gamma_\eta$
is bounded above by\break 
\hbox{$\mu_\xo(\rand)/\mu_\xo(\tilde\eta)$}.
Since $\mu_\xo$
is a finite measure and $\mu_\xo(\tilde\eta)>0$, the above sum
converges.\qed

\begin{thr}\label{radpoint}
If $\delta_\theta(\Gamma)>0$ then a regular radial limit point $\tilde\eta\in \radlim\cap \regrand$ is not a
point mass for any \bd y. 
\end{thr}
\prf\  Let $\mu$ be a \bd y. If $\tilde\eta\notin
\rand_\theta$, then $\tilde\eta\notin\supp (\mu_\xo)$, hence $\tilde\eta$ cannot be a point
mass.

Suppose $\tilde\eta=(\eta_1,\eta_2,\theta)\in\radlim\cap \rand_\theta$ is a point mass for $\mu$. Then by Theorem~\ref{posmeasrad}  
$\langle b,H_\theta\rangle=\delta_\theta(\Gamma)>0$, hence by continuity of the map $\hat\theta\mapsto \langle b,H_{\hat\theta}\rangle$ there exists $\eps>0$ \st  every $\hat\theta\in [0,\pi/2]$ with $|\hat\theta-\theta|<\eps$ satisfies $\langle b,H_{\hat\theta}\rangle >q>0$. Moreover, by definition of the radial limit set~(\ref{radlimthetacover})  there exists  a constant
$c>0$ and a sequence $(\gamma_n)=\big((\gamma_{n,1}, \gamma_{n,2})\big)\subset\Gamma$  \st $|\theta(\xo,\gamma_n\xo)-\theta| <\eps$ and $\tilde\eta\in S(o:
B_{\gamma_n o}(c))$ for all $n\in\NN$.
Corollary \ref{esti}  implies
$\bs_{\eta_i}(o_i,\gamma_{n,i} o_i)>d_i(o_i,\gamma_{n,i} o_i)-2c$ for all $n\in\NN$ and
$i\in\{1,2\}$. 
We  conclude
$$ \bs^b_{\tilde\eta}(\xo,\gamma_n \xo)>\langle b, H(o,\gamma_n o) \rangle-2 \Vert b\Vert_1 c  \to \infty\,,$$
because $\langle b, H(\xo,\gamma_n\xo)\rangle >q\cdot d(\xo,\gamma_n\xo)$ and $d_i(o_i,\gamma_{n,i} o_i)\to\infty$ for all $i\in\{1,2\}$ as
$n\to\infty$. Passing to a subsequence if necessary, we may therefore
assume that 
$\ \bs^b_{\tilde\eta}(o,\gamma_{n} o)\rangle\  $
is strictly increasing to infinity as $n\to\infty$.

Now suppose there exist
$l,j \in\NN$, $l\ne j$ \st 
$\gamma_l^{-1}\Gamma_{\tilde\eta} = \gamma_j^{-1}\Gamma_{\tilde\eta}$.
Since $\tilde\eta$ is a point mass for $\mu$  Lemma~\ref{busnull} implies
$$\bs^b_{\tilde\eta}(o,\gamma_j o) =\bs^b_{\tilde\eta}(o,\gamma_{l}\xo)\,,$$
in contradiction to the choice of the subsequence $(\gamma_{n})$. 
Hence
\hbox{$\gamma_l^{-1}\Gamma_{\tilde\eta}{\neq}\gamma_j^{-1}\Gamma_{\tilde
\eta}$} for all $l\ne j$, and
the sum $\sum e^{\bs^b_{\tilde\eta}(\xo,\gamo)}$ over a system of coset representatives of $\Gamma/\Gamma_{\tilde\eta}$ is
bounded below by 
$$\sum_{n\in\NN}e^{\bs^b_{\tilde\eta}(o,\gamma_{n}\xo)}, $$
and therefore diverges in contradiction to {\rm Lemma}
\ref{sumconv}. We conclude
that  $\tilde\eta$ cannot be a point mass for
$\mu_\xo$.\qed

\section{Hausdorff dimension}\label{Hausdorff}
% For the remainder of the article  $\XX$ is a product of locally compact %geodesically complete 
% Hadamard spaces $\XX_1$, $\XX_2$, $\xo=(\xo_1, \xo_2)$ a fixed base point, and $\Gamma\subset\is(\XX_1)\times\is(\XX_2)$ a discrete group which
%  contains two isometries $g=(g_1,g_2)$ and $h=(h_1,h_2)$ \st $g_i$ and $h_i$ are independent rank one elements of $\Gamma_i$  for $i=1,2$.  
% 
This final section introduces an appropriate notion of Hausdorff measure and Hausdorff dimension on the geometric boundary $\rand$ in order to estimate the size of the radial limit set in each $\Gamma$-invariant subset $\Lim\cap\rand_\theta$ of the  geometric limit set. Our results are most precise for a class of groups which we  call radially cocompact. In this case, the Hausdorff dimension of the radial limit set in a
given subset $\rand_\theta\subseteq\regrand$ equals the exponent of growth of slope $\theta$

We will follow the idea of G.\ Knieper (\cite[\S4]{MR1465601})
for a
definition of Hausdorff measure on the geometric boundary.
For $\tilde\xi\in\rand$, $c>0$ and $0< r<e^{-c}$ we call the set
$$ 
B_r^c(\tilde\xi):=\big\{\tilde\eta\in\rand : d(\sigma_{\xo,\tilde\eta}(-
\log
r),\sigma_{\xo,\tilde\xi}(-\log r))<c\big\}$$ 
a {\hl $c$-ball} of radius $r$ centered at $\tilde\xi$. 
Using this conformal structure, we define as in the case of metric spaces
Hausdorff measure and Hausdorff dimension on the geometric boundary $\rand$.

\begin{df}
Let $E$ be a Borel subset of $\rand$, and
$$ 
\Hd_\varepsilon^{\alpha}(E):=\inf\Big\{\sum r_i^\alpha : 
|E\subseteq\bigcup B_{r_j}^c(\tilde\xi_j)\,,\ r_j<\varepsilon
\Big\}\,.$$
The  $\alpha$-dimensional {\hd Hausdorff measure} of $E$ is defined by  
\mbox{$\Hd^{\alpha}(E)=\displaystyle\lim_{\varepsilon\to 0}\Hd_
\varepsilon^{\alpha}(E)$}, and 
the {\hd Hausdorff dimension} of $E$ is the number
$$
\Hdim(E)=\inf\big\{\alpha \ge 0\bigm|\Hd^{\alpha}(E)<
\infty\}\,.$$
\end{df}

Recall the notion of Weyl chambers and Weyl chamber shadows from Section~\ref{prodHadspaces}. 
The following lemma gives  a 
relation between Weyl chamber shadows in $\rand_\theta$ and $c$-balls.  

\begin{lem}\label{shadowinball}
Let $c>0$ and $\theta\in [0,\pi/2]$.  We set
$A:=\max\{\sin\theta/\cos\theta, \cos\theta/\sin\theta\}$ if $\theta\in (0,\pi/2)$ and $A:=0$ otherwise. 
%Then there exists $R_0>c$ \st the following holds:
Let  $\eps\in (0,\pi/2)$ arbitrary with $\eps\le \frac12 \min\{\theta, \pi/2-\theta\}$ if $\theta\in (0,\pi/2)$, and $\tilde\eta\in\rand_\theta$.
If  $y\in \Ch_{\xo,\tilde\eta}$ satisfies $d(\xo,y)>c$ and $|\theta(\xo,y)-\theta|<\eps$, then with $r:=\exp(-d(\xo,y)(\cos \eps -A\sin\eps))$ we have
$$
\mathrm{Sh}\big(\xo:B_y(c/2)\big)\cap \rand_\theta\subseteq B^c_r
(\tilde\eta)\,.$$
\end{lem}
\prf\  Fix $y=(y_1,y_2)\in \Ch_{\xo,\tilde\eta}$ with $d(\xo,y)>c$ and $|\theta(\xo,y)-\theta|<\eps$, and set $t_i:=d_i(\xo_i,y_i)$ for $i\in\{1,2\}$. We have to show that for 
$\tilde\zeta\in \Sh\big(\xo:B_y(c/2)\big)\cap \rand_\theta$ arbitrary the inequality $\ d(\sigma_{\xo,\tilde\eta}(-\log r), \sigma_{\xo,\tilde\zeta}(-\log r))<c\,$ holds. Notice
that by 
definition of the Weyl chamber shadows 
we have $d(y,\Ch_{\xo,\tilde\zeta}) <c/2$. 

Assume first that $\theta\in (0,\pi/2)$ and write $\tilde \eta=(\eta_1,\eta_2,\theta)$, $\tilde\zeta=(\zeta_1,\zeta_2,\theta)$. 
If for $i\in\{1,2\}$ we set $c_i:=d_i(y_i,\sigma_{\xo_i,\zeta_i})$, then $d(y,\Ch_{\xo,\tilde\zeta})<c/2$ implies $\sqrt{c_1^2+c_2^2}<c/2$. 
Since $y\in \Ch_{\xo,\tilde\eta}$  $\ y_1$ is a point on the geodesic ray $\sigma_{\xo_1,\eta_1}$ and $y_2$ is a point on the geodesic ray $\sigma_{\xo_2,\eta_2}$. 
Moreover, by elementary geometric estimates  we have  for $i=1,2$
$$d_i(\sigma_{\xo_i,\eta_i}(t_i), \sigma_{\xo_i,\zeta_i}(t_i))=d_i(y_i,\sigma_{\xo_i,\zeta_i}(t_i)) \le 2c_i\,.$$
Now $|\theta(\xo,y)-\theta|<\eps$ and the choice of $A$ imply 
\be t_1&= & d(\xo,y)\cos\theta(\xo,y)\ge d(\xo,y)\cos(\theta+\eps)\ge d(\xo,y)\cos\theta\big(\cos\eps-A\cdot \sin\eps\big)\,,\\
t_2&= & d(\xo,y)\sin\theta(\xo,y)\ge d(\xo,y)\sin(\theta-\eps)\ge d(\xo,y)\sin\theta\big(\cos\eps -A\cdot\sin\eps\big)\,.\ee
Notice  that $\cos\eps -A\cdot\sin\eps >0$ because $\eps\le \frac12 \min\{\theta,\pi/2-\theta\}$. Using the definition of the constant  $r$ and the convexity of the distance function in $\XX_1$, $\XX_2$ 
we get
\be d_1(\sigma_{\xo_1,\eta_1}(-\log r \cdot \cos\theta), \sigma_{\xo_1,\zeta_1}(-\log r \cdot \cos\theta)) & \le & d_1(\sigma_{\xo_1,\eta_1}(t_1), \sigma_{\xo_1,\zeta_1}(t_1))\le 2 c_1\,,\\
d_2(\sigma_{\xo_2,\eta_2}(-\log r \cdot \sin\theta), \sigma_{\xo_2,\zeta_2}(-\log r \cdot \sin\theta)) & \le & d_2(\sigma_{\xo_2,\eta_1}(t_2), \sigma_{\xo_2,\zeta_2}(t_2))\le 2 c_2\,.\ee
Since $\sigma_{\xo,\tilde\eta}(t)=\big(\sigma_{\xo_1,\eta_1}(t\cos\theta),\sigma_{\xo_2,\eta_2}(t\sin\theta)\big)$ and  
$\sigma_{\xo,\tilde\zeta}(t)=\big(\sigma_{\xo_1,\zeta_1}(t\cos\theta),\sigma_{\xo_2,\zeta_2}(t\sin\theta)\big)$ for all $t>0$, we conclude
$$ d(\sigma_{\xo,\tilde\eta}(-\log r), \sigma_{\xo,\tilde\zeta}(-\log r)) \le \sqrt{(2c_1)^2+(2c_2)^2}=2\sqrt{c_1^2+c_2^2}\le c\,.$$

If $\theta=0$ we write $\tilde \eta=\eta_1$ and $\tilde\zeta=\zeta_1$. By (\ref{shadowsingular}) we have $d_1(y_1,\sigma_{\xo_1,\zeta_1})<c/2$, hence by the same reasoning as before 
$$d_1(\sigma_{\xo_1,\eta_1}(t_1), \sigma_{\xo_1,\zeta_1}(t_1))=d_1(y_1,\sigma_{\xo_1,\zeta_1}(t_1)) \le c\,.$$
Moreover,  $t_1=  d(\xo,y)\cos\theta(\xo,y)\ge d(\xo,y)\cos\eps= d(\xo,y)\big(\cos\eps-A\cdot \sin\eps\big)$. The conclusion then follows from the convexity of the distance function as above.

The case $\theta=\pi/2$ is analogous; we only have to notice that \vspace*{-2mm}
 $$t_2=d(\xo,y)\sin\theta(\xo,y)\ge d(\xo,y)\sin(\pi/2-\eps)=d(\xo,y)\cos\eps= d(\xo,y)\big(\cos\eps-A\cdot \sin\eps\big).$$
\qed\\

The inclusions from the previous lemma allow to give an upper bound for the Hausdorff dimension
of the radial limit set. 
\begin{thr}\label{upbound}
If $\theta\in [0,\pi/2]$, then the Hausdorff dimension of the radial
limit set in $\rand_\theta$ 
is bounded above by $\delta_{\theta}(\Gamma)$.
\end{thr}
\prf\    Let $\theta\in [0,\pi/2]$ and fix  $c>0$ sufficiently large. By definition of the radial limit set
$$\radlim\cap \rand_\theta\subseteq \bigcap_{R>c}\, \bigcap_{\eps>0}
\bigcup_{\begin{smallmatrix}{\scriptscriptstyle\gamma\in
\Gamma}\\{\scriptscriptstyle d(o,\gamma\xo)>R}\\
{\scriptscriptstyle |\theta(\xo,\gamma\xo)-\theta|<\eps}
\end{smallmatrix}}
\Sh(\xo: B_{\gamma\xo}(c/2))\,. $$
Fix $A\ge 0$ as in Lemma~\ref{shadowinball} and let $\eps\in (0,\pi/2)$ arbitrary, $\eps<\frac12\min\{\theta,\pi/2-\theta\}$ if $\theta\in (0,\pi/2)$. Put 
$\hat\Gamma:=\{\gamma\in\Gamma : \Sh(\xo:B_{\gamma\xo}(c/2))\cap\radlim\cap\rand_\theta\ne\emptyset,\, |\theta(\xo,\gamma\xo)-\theta|<\eps\}$.
For  $\gamma\in\hat\Gamma$  we denote $\tilde\xi_\gamma$ a point in $\Sh(\xo:B_{\gamma\xo}(c/2))\cap\radlim\cap\rand_\theta$ and set
$$ r_\gamma:=\exp\big(-d(\xo,\gamma\xo)(\cos\eps-A \sin\eps)\big)\,.$$
Let $\rho<e^{-c}$ and set $\Gamma':=\{\gamma\in\hat\Gamma : r_\gamma<\rho\}$. By the previous lemma we have 
$\Sh(\xo:B_{\gamma\xo}(c/2))\cap \rand_\theta\subseteq B^
c_{r_\gamma}(\tilde\xi_{\gamma})$ 
for all $\gamma\in\Gamma'$, 
hence 
$$\radlim\cap \rand_\theta\subseteq \bigcup_{\gamma\in
\Gamma'}\
   B^c_{r_\gamma}(\tilde\xi_{\gamma})\,.$$
Using the definition of $\Hd_\rho^{\,\alpha}$ we estimate
$$ \Hd^{\alpha}_\rho(\radlim\cap \rand_\theta)\le\sum_{\gamma\in
\Gamma'}r_\gamma^\alpha=\sum_{\gamma\in\Gamma'} e^{-\alpha
(\cos\eps-A  \sin\eps)d(\xo,\gamma\xo)}\,.$$
Recall from Section~\ref{ExpGrowth} that 
$$
Q_\theta^{s,\eps}(\xo,\xo):= \sum_{\begin
{smallmatrix}{\gamma\in\Gamma}\\{|\theta(\xo,\gamma\xo)-\theta|<\eps}\end{smallmatrix}} e^{-s d(o,\gamma
  \xo)} $$ converges for $s >\delta_{\theta}^{\eps}(o,\xo)$. Hence if
  $\hat s:=
  \alpha(\cos\eps-A\sin\eps)>\delta_{\theta}^{\eps}(o,\xo)$,
  we have 
$$\Hd^{\alpha}_\rho(\radlim\cap \rand_\theta)\le 
Q_{\theta}^{\hat s,\eps}(\xo,\xo)<\infty\,.$$
This shows that for
$\alpha>\delta_{\theta}^{\eps}(o,\xo)/(\cos\eps-A\sin\eps)$,
$\Hd^{\alpha}_\rho(\radlim\cap \rand_\theta)$ is finite. Taking the limit as
$\eps\searrow 0$, we conclude that the same is true for $\alpha>
\delta_{\theta}(\Gamma)$. Letting $\rho\searrow 0$, we obtain
$\, \Hd^{\alpha}(\radlim\cap \rand_\theta)<\infty\,\ $  if $\alpha>
 \delta_{\theta}(\Gamma),\ $ hence $\ \Hdim(\radlim\cap\rand_
\theta)\le \delta_{\theta}(\Gamma)$.\qed\\[-1mm]

Notice that in the previous proof we only used the definition of $\delta_\theta(\Gamma)$ and not the existence of a $(b,\theta)$-density. In particular, the claim also holds for
slopes $\theta\in [0,\pi/2]$ for which $\delta_\theta(\Gamma)=0$. 

The notion of convex cocompact and geometrically finite groups plays
an important role in the theory of Kleinian groups. For these groups 
the Hausdorff dimension of the limit set is equal to the critical
exponent. We suggest here the
following definition to replace convex cocompactness in a more general setting.
\begin{df}\label{radcoc} If $\XX$ is a locally compact Hadamard space, then a  discrete group $\Gamma\subset\is(\XX)$ is called
  {\hd radially cocompact} if there exists a constant $c_\Gamma>0$ \st
  for any $\tilde\eta\in\radlim$ and for all $t>0$ there exists an element
  $\gamma\in\Gamma$ with 
$$ d\big(\gamma\xo,\sigma_{\xo,\tilde\eta}(t)\big)<c_\Gamma\,.$$
\end{df}
The most familiar radially cocompact groups are convex cocompact
isometry groups of rank one symmetric spaces and  uniform
lattices of higher rank symmetric spaces or Euclidean buildings. A further example is given 
by products of convex cocompact groups acting on the Riemannian
product of two Hadamard spaces with pinched negative curvature. 

For radially cocompact discrete groups $\Gamma\subset \is(\XX_1)\times\is(\XX_2)$, the existence of a \bd y $\mu$
together with {\rm Theorem} \ref{shadowlemma} allows to obtain a
lower bound for the Hausdorff dimension of the radial limit set in  $\rand_\theta$. From
here on, we fix $c>2\max\{c_\Gamma,c_0\}$ with $c_\Gamma$
as in Definition \ref{radcoc} and $c_0$ as in Theorem
\ref{shadowlemma}. 

\begin{thr}\label{lowbound}
Let $\Gamma\subset\is(\XX)$ be radially cocompact, $\theta\in [0,\pi/2]$, and  $\mu$ a \bd y. Then 
there exists a constant $C_0>0$ \st for any Borel subset $E\subseteq \radlim$
$$ \Hd^{\langle
b,H_\theta\rangle}(E)\ge C_0\cdot \mu_\xo(E)\,.$$
\end{thr}
\prf\  Set $\alpha:=\langle b, H_\theta\rangle$. Since $\Hd^{\alpha}(E)\ge \Hd^{\alpha}(E\cap \rand_\theta)$ and
$\mu_\xo(E)=\mu_\xo(E\cap \rand_\theta)$, it suffices to prove the assertion for
$E\subseteq\radlim\cap \rand_\theta$. Let $\rho>0$, $q>0$ arbitrary,
and choose a cover of $E$ by balls
$B_{r_n}^{c}(\tilde\eta_n)$, $r_n<\rho$,  \st 
$$\Hd_\rho^{\alpha}(E)\ge\sum_{n\in\NN}r_n^{\alpha}-q\,.$$
If $ B_{r_n}^{c}(\tilde\eta_n)\cap E=\emptyset$, we do not need
$B_{r_n}^{c}(\tilde\eta_n)$ to cover $E\subseteq\radlim\cap \rand_\theta$, otherwise we choose $\tilde\xi_n\in
B_{r_n}^{c}(\tilde\eta_n)\cap E$. Since
$\Gamma$ is radially cocompact, there exists $\gamma_n=(\gamma_{n,1},\gamma_{n,2})\in\Gamma$ \st
\begin{equation}\label{abstandradcoc}
d(\gamma_n\xo, \sigma_{\xo,\tilde\xi_n}(-\log r_n))\le c\,. 
\end{equation}
If $\theta\in (0,\pi/2)$ we write $\tilde\xi_n=(\xi_{n,1},\xi_{n,2},\theta)$, if $\theta=0$ we set $\xi_{n,1}=\tilde\xi_n\in\rand_1$ 
and choose $\xi_{n,2}\in\rand_2$ arbitrary, if $\theta=\pi/2$ we let $\xi_{n,1}\in\rand_1$ arbitrary and set $\xi_{n,2}=\tilde\xi_n\in\rand_2$.
With this notation  inequality~(\ref{abstandradcoc}) implies
$$ d_1(\gamma_{n,1}\xo_1,\sigma_{\xo_1,\xi_{n,1}}(-\log r_n \cdot \cos\theta))\le c\quad\mbox{and}\qquad d_2(\gamma_{n,2}\xo_2,\sigma_{\xo_2,\xi_{n,2}}(-\log r_n \cdot \sin\theta))\le c\,,$$
so using the triangle inequalities we obtain
\be 
\big| d_1(\xo_1,\gamma_{n,1}\xo_1)-d_1\big(\xo_1, \sigma_{\xo_1,\xi_{n,1}}(-\log r_n \cdot \cos\theta)\big)\big|&\le &c \qquad\mbox{and}\\
\big| d_2(\xo_2,\gamma_{n,2}\xo_2)-d_2\big(\xo_2, \sigma_{\xo_2,\xi_{n,2}}(-\log r_n \cdot \sin\theta)\big)\big|&\le &c\,.
\ee
We therefore estimate
\be \langle b,H(\xo,\gamma_n\xo)\rangle &\ge& -|b_1|c+  b_1 d_1\big(\xo_1, \sigma_{\xo_1,\xi_{n,1}}(-\log r_n \cdot \cos\theta)\big)  -|b_2|c\\
&& + b_2 d_2\big(\xo_2, \sigma_{\xo_2,\xi_{n,2}}(-\log r_n \cdot \sin\theta)\big)= -\langle b,H_\theta\rangle \log r_n - \Vert b\Vert_1 c,
\ee
hence
$$-\langle b, H(\xo,\gamma_n\xo)\rangle \le \langle b, H_\theta\rangle\log r_n +c\Vert b\Vert_1  \,.$$
Furthermore, we have  $B_{r_n}^{c}(\tilde\eta_n)\subseteq \Sh(\xo: B_{\gamma_n\xo}(3c))$,  % for any $\rho>0$,
hence
\hbox{$E\subseteq \bigcup_{n\in\NN} \Sh(\xo: B_{\gamma_n\xo}(3c))$}.
We conclude 
\begin{align*}
\mu_\xo(E)&\le\mu_\xo\Big(\bigcup_{n\in\NN}\Sh\big(\xo: B_{\gamma_n\xo}(3c)\big)\Big)
\le \sum_{n\in\NN}  \mu_\xo\big(\Sh(\xo:
B_{\gamma_n\xo}(3c))\big)\\
&\le D(3c) \sum_{n\in\NN} e^{-\langle b, H(\xo,\gamma_n\xo\rangle)}
 \le D(3c) \sum_{n\in\NN} e^{\alpha\log r_n +c\Vert b\Vert_1 }\\
&\le D(3c) e^{c \Vert b\Vert_1} \sum_{n\in\NN}
 r_n^{\alpha} \le  D(3c) e^{c\Vert b\Vert_1} \big(
   \Hd_\rho^{\,\alpha}(E)+q\big)\,.
\end{align*} 
The claim now follows as $q\searrow 0$ and $\rho\searrow 0$.\qed\\[-1mm]
\begin{thr}
Let $\Gamma\subset\is(\XX_1)\times \is(\XX_2)$ be  radially cocompact, and $\theta\in (0,\pi/2)$ \st $\delta_\theta(\Gamma)>0$. Then $$\Hdim(\radlim\cap\rand_\theta) = \delta_{\theta}(\Gamma)\,.$$ 
\end{thr} 
\prf\  From Section~\ref{genPseries} we know that there exists a $(b,\theta)$-density $\mu$ for appropriate parameters $b=(b_1,b_2)$. So the previous theorem implies that for 
$\alpha:= \langle b,  H_\theta \rangle$ 
 $$\Hd^{\alpha}(\radlim\cap \rand_\theta)\ge  C_0 \mu_\xo(\radlim)\ge 0,$$
hence
$$\Hdim(\radlim\cap\rand_\theta)\ge\alpha=\langle b,H_\theta\rangle\ge \delta_{\theta}(\Gamma)$$ 
%If $\theta\in (0,\pi/2)$, then 
 by Theorem~\ref{deltaklein}. The assertion now follows directly
from Theorem~\ref{upbound}.\qed\\[-1mm]

For Example 1 described in Section~\ref{ExpGrowth} we deduce
the following 
\begin{cor}
Let $\XX=\XX_1\times \XX_2$ be a product of two Hadamard manifolds of pinched negative curvature, $\Gamma_1\subset \is(\XX_1)$, $\Gamma_2\subset \is(\XX_2)$  convex cocompact groups with
critical exponents $\delta_1$, $\delta_2$, and $\Gamma=\Gamma_1\times \Gamma_2\subset\is(\XX)$. Then for any $\theta\in (0,\pi/2)$
we have 
$$\Hdim(\radlim\cap
\rand_\theta) = \delta_1\cos\theta+\delta_2\sin\theta\,.$$
\end{cor}

\bibliography{Bibliographie} 

\vspace{1.5cm}
\noindent Gabriele Link\\
Institut f\"ur Algebra und Geometrie\\
Karlsruher Institut f\"ur Technologie (KIT)\\
Kaiserstr. 89-93  \\[1mm]
D-76133 Karlsruhe\\
e-mail:\ gabriele.link@kit.edu

\end{document}